# On Jordan's measurements

**Frédéric Brechenmacher**

Univ. Lille Nord de France, F-59000 Lille

U. Artois. Laboratoire de mathématiques de Lens (EA 2462),

rue Jean Souvraz, S.P. 18, F-62300 Lens, France

frederic.brechenmacher@euler.univ-artois.fr

## Introduction

The Jordan measure, the Jordan curve theorem, as well as the other generic references to Camille Jordan's (1838-1922) achievements highlight that the latter can hardly be reduced to the "great algebraist" whose masterpiece, the *Traité des substitutions et des equations algébriques*, unfolded the group-theoretical content of Évariste Galois's work. Not only did Jordan also write the influential *Cours d'analyse de l'École polytechnique* (1882-1887) [Gispert 1982], but more than two-third of his publications were not classified in algebra by his contemporaries.

The present paper appeals to the database of the reviews of the *Jahrbuch über die Fortschritte der Mathematik* (1868-1942) for providing an overview of Jordan's works. On the one hand, we shall especially investigate the collective dimensions in which Jordan himself inscribed his works (1860-1922). On the other hand, we shall address the issue of the collectives in which Jordan's works have circulated (1860-1940).

At the turn of the 20th century, Jordan has been assigned a specific role in regard with Galois's legacy.[1] When he died in 1922, this relation to Galois was at the core of Jordan's identity as an algebraist.[2] Later on, in the second half of the 20th century, several authors pointed out that the relation Jordan-Galois did not fit the historical developments of group theory in connection to Galois theory. As a result, and despite the key positions Jordan occupied in the Parisian mathematical scene, the latter has thus often been considered as isolated, and in fact almost foreign to the French mathematical scene.[3] Moreover, the reception of the *Traité* has been believed to be very limited.[4]

Some recent works have nevertheless shown than rather than Jordan's position or Jordan's book, it is the usual presentation of the Jordan-Galois relation that does not fit History [Brechenmacher 2011a]. The traditional role assigned to Jordan indeed results from both the retrospective perspective of group theory and from the public dimension of the figure of Galois [Brechenmacher 2012]. These new perspectives call for a new look on the collective dimensions of Jordan's works.

---

[1] Cf. [Lie 1895]; [Picard 1897].
[2] Cf. [Adhémar 1922]; [Picard 1922]; [Lebesgue 1923]; [Villat 1922]; [Taton 1947, p.94].
[3] Cf. [Klein 1921, p. 51,] [Julia *in* Galois 1962], [Dieudonné 1962].
[4] Cf. [Wussing 1984], [Kiernan 1971].

Moreover, the time-period during which Jordan has been publishing his works, i.e., 1860-1922, provides an opportunity to investigate some collective organizations of knowledge that pre-existed the development of object-oriented disciplines such as group theory (Jordan-Hölder theorem), linear algebra (Jordan's canonical form), topology (Jordan's curve), integral theory (Jordan's measure), etc.

At the time when Jordan was defending his thesis in 1860, it was common to appeal to transversal organizations of knowledge, such as what the latter designated as the "theory of order." When Jordan died in 1922, it was however more and more common to point to object-oriented disciplines as identifying both a corpus of specialized knowledge and the institutionalized practices of transmissions of a group of professional specialists. Jordan, indeed, was one of the first French academic mathematicians to be celebrated as an "algebraist."[5]

It is the aim of the present paper to investigate in parallel both Jordan's works and the corpus of the *Jahrbuch* reviews that referred to "Jordan" between 1868 and 1939 (the Jordan corpus, for short). It should be pointed out at this stage that such a corpus cannot be automatically generated. Indeed one has to remove manually from the database the reviews that are referring to homonyms (Wilhelm Jordan, Paul Jordan, Charles Jordan etc.) as well as the ones that refer to the retrospective sporadic additions of the classification AMS 2000 (Jordan's curve, Jordan's canonical form, etc.).

Let us now reflect on the appropriate measurements for the Jordan corpus. Two time-periods have to be distinguished.

First, the number of reviews of the *Jahrbuch* involved from 1868 to 1910 (i.e., 275 reviews referring to Jordan plus about 100 reviews on papers published by Jordan himself) makes it possible to investigate each review individually. Here, the role of the diagrams given in the annex to this paper will be to support our presentation of the Jordan corpus. Moreover, the *Jahrbuch* database will be used for approximating the solution to the inverse intertextual relations problem.[6] Unlike the unknown corpus of all the papers that referred to Jordan from 1868 to 1910, the corpus of the *Jahrbuch* reviews can be investigated systematically. It is especially possible to appeal to the intertextual relations indicated by these reviews for investigating how the texts under reviews were fitting some collective networks of texts.[7]

To be sure, the present paper is nevertheless not directly investigating the circulations of Jordan's works but the echoes of these circulations in the *Jahrbuch*. It may indeed happen that one reviewer may add some references of his own, or that one may omit to mention some of the main references of the paper under review. As a matter of fact, some papers that are pointing to Jordan's works as one of their main references do not appear in the Jordan corpus.[8] Some of these were actually not even reviewed in the *Jahrbuch*.[9] For the purpose of tracing down such situations, we shall consider a few sub-corpora as probes. In connection to these sub-corpora, the intertextual relations of the *papers* under reviews will be systematically investigated in parallel to the ones of the *reviews*.

We will also investigate other collective dimensions such as the distributions of the section of the *Jahrbuch* mathematical classification. As has been shown in [Goldstein, 1999, p. 204-212], such distributions usually do not fit the intertextual organizations of texts into networks. The sections of the mathematical classification indeed usually encompass various types of works. They nevertheless reflect some collective opinions on the organization of mathematics in sub-domains and therefore deserve to be studied as such.[10] In what follows, I shall appeal to some abbreviated denominations of the chapters of the *Jahrbuch*. The full names will be indicated when these chapters will be mentioned for the first time. For instance, until the mid 1890s the third chapter of the section "II. Algebra" mixed together various issues related to substitutions, determinants and invariants. For short, I will designate this chapter as the "substitutions" chapter. One shall nevertheless keep in mind that this chapter encompassed a much broader spectrum that the theory of substitutions. An list of the *Jahrbuch* chapter involved here is provided in the annex to this paper.

The second time-period under consideration, i.e., the period from 1910 to 1940, involves more than 700 texts, thereby making it possible to appeal to some quantitative analysis. Moreover, unlike the publications they point to, the *Jahrbuch* reviews provide a corpus homogeneous enough for the use of quantitative methods. We shall thus focus on counting the reviews in the Jordan corpus in contrast with some other quantitative data that would be more specific to the papers under review (number of pages, number of journals involved, number of authors etc.). It must nevertheless be pointed out that the changing and heterogeneous practices of both publications and references make it customary to appeal to some qualitative analysis of the corpus in addition to the quantitative ones.

To be sure, our investigations of the post-1910 corpus shall stay at a more superficial level than our analysis of the pre-1910 period. But in any case, in addition to a quantitative shift in the Jordan corpus, the status of the references to Jordan also changed after 1910. Generic references indeed play a more and more important role. To be sure, generic designations raise different issues than direct references. It is indeed usually pointless to wonder about which specific part of which specific paper in Jordan's works, a reference to the Jordan curve theorem is pointing to. Generic designations rather raise some issues on the evolutions of the concept, the method, or the theorem they are pointing to. Although it is not the place here to investigate such issues, we shall nevertheless present a case study on the process of universalization of the generic designation "Jordan's canonical form theorem."

---

[8] Such as the paper in which Hamburger (himself a *Jahrbuch* reviewer) pointed out in 1872 the connections between Jordan's canonical form and Weierstrass' elementary divisors theorem [Brechenmacher 2007].
[9] For instance, the famous paper in which Hölder formulated the Jordan-Hölder theorem by appealing to the notion of quotient group was not reviewed in the *Jahrbuch* [Hölder 1889]
[10] See [Goldstein 1999] and [Siegmund Schulze 1993].

Finally, given the number of (reviews on) publications considered in the present paper, the delimitation of a relevant bibliography raises a specific difficulty. Indeed, any single paper among the hundreds of texts involved in the Jordan corpus deserves to be identified precisely. To be sure, most of these publications shall nevertheless not be listed in the present paper. Even though the bibliography will be limited to a few representative texts, I shall however refer to a greater number of papers by providing enough information to identify these in the *Jahrbuch* database of the Jordan corpus (usually in the form (Name Year)).

## 1. An overview on Jordan's mathematical works

The present overview mainly aims at introducing the analysis of the various collective forms of references to Jordan's works that it is developed in the next sections of this paper.[11] This first section therefore retrospectively focuses on the same aspects of Jordan's works as the *Jahrbuch* reviews of the Jordan corpus. But I shall nevertheless aim at shedding some new light on some collective organizations of texts that lie beneath the global corpus of Jordan's works without appealing to retrospective categories such as group theory, abstract algebra, topology, etc.

Jordan published most of his research papers before the 1880s. On the other hand, some academic, teaching, and editorial responsibilities dominate most post-1881 publications. The corpus of Jordan's works will therefore be divided into the following two groups of time-periods: ([1860, 1867], [1868, 1874], [1875, 1880]) and ([1881, 1904], [1905, 1922]).

The diagrams n°1 represent the distribution of the classifications of Jordan's papers for each post 1868 time-period (recall that the *Jahrbuch* was founded in 1868). The chapters "II. Algebra / 3. Substitutions" and "II. Algebra / 1.Equations" dominated the time-period 1868-1874. From 1875 to 1880, the evolution was mainly due to the increasing weights of the chapters "II. Algebra / 2.Theory of forms" and "VI. Differential and integral calculus / 5.Differential equations." The various editions of the three volumes of the *Cours d'analyse* dominated the next twenty-four-year time-period. After he had left in 1912 his teaching duties at both the *École polytechnique* and the *Collège de France* (to Jacques Hadamard and Georges Humbert respectively), Jordan published again some papers on various topics, such as on the investigation of the algebraic forms invariant by some given groups of substitutions.

Let us now consider more closely each of these time-periods.

---

[11] For some more mathematical details on Jordan's works, see [Dieudonné 1962].

### 1.1. The theory of order (1860-1867) - the regular solids research field

Most of Jordan's works from 1860 to 1867 fit into what the author designated as "the theory of order" in his 1860 thesis, and again in his 1881 application to the *Académie* [Jordan, 1881, p. 7-8].

This theory was characterized as the part of mathematics dealing with "relations" between classes of objects in contrast with classical concerns for quantities, magnitudes, or proportions. It was thus transversal to algebra, number theory, geometry and kinematics. Jordan especially appealed to the long-term legacy of Louis Poinsot. From 1808 to 1844, the latter had highlighted several times the transversal role played by the notion of "order" in the analogies encountered in various cyclic situations such as the investigations of cyclotomic equations, congruences, symmetries, polyhedrons, and motions. Poinsot had especially discussed the notion of order when commenting on Gauss's number-theoretic indexing of the roots of the cyclotomic equations of the division of the circle. Later on, he also characterized the theory of order as having a relation to algebra analogous to the relations between Gauss's higher arithmetic and usual arithmetic, or that between *analysis situs* and geometry [Boucard 2011].

The theory of order encompasses the variety of issues Jordan tackled from 1860 to 1867: substitutions groups, algebraic equations, higher congruences, motions of solid bodies, symmetries of polyhedrons, crystallography, the *analysis situs* of deformations of surfaces (including Riemann surfaces), and the groups of monodromy of linear differential equation to which Jordan's second thesis was devoted. As an example of the intertwining of such various topics, one may cite Jordan's 1868 memoir on the classification of the groups of motions of solid bodies in the legacy of Bravais' crystallography.

As shall be seen later, this specific mix of topics had underlying it the larger-scale collective dimension of a field of research at a European level. This organization was nevertheless not commonly identified as the theory of order and I shall designate it by the *regular solid research field*, thereby pointing to a pattern most of the authors of the field appealed to. To be sure this field should be investigated further. Ehrhard Scholz's approach to the history of crystallography sheds light on some important aspects of the *regular solid research fields*, such as the French approach to crystallography through the investigations of both the symmetries and motions of polyhedrons as well as through some topological considerations on the deformations of networks of points in connection the arithmetic of quadratic forms [Scholz 1989, p. 81-110]. But some other aspects require further investigations, such as the connections of Jordan's theory of order with some approaches to mechanics at the *École polytechnique* or with the long-terme tradition of the cinematic approach to geometry.[12]

---

[12] Méray's 1874 unified cinematic approach to both space and plane geometry [Bkouche 1991] seems especially close to Jordan's contemporary $n^{th}$ dimensional agenda. Moreover, Houël's 1867 claims on the importance to distinguish between the abstract notion of motion as it may be used in mechanics seems also to resound in some of Jordan's claims on groups. Finally, the way Chasles advocated the abstraction and the generality of synthetic geometry when he investigated the motions of solid bodies is close to Jordan's topological considerations. Chasles's works would be referred to later on by Schönflies in connection to the Jordan curve theorem.

### 1.2. The "method of reduction" and the analytic representation of linear substitutions

In his first thesis, Jordan had already sown the seeds of the group theoretical organization of knowledge by which he would replace the theory of order in the late 1860s. The main result of the thesis was the introduction of the *group of n-ary linear substitutions on integers mod.p* (i.e. $Gl_n(F_p)$ in nowaday's parlance). As shall be seen later, this introduction would give rise to what would be one of the main ways to refer to Jordan at the turn of the century, i.e., to "Jordan's linear groups in Galois fields."

We shall thus take a closer look at this episode. Jordan's first thesis was devoted to the problem of the number of values of functions under the action of substitutions.[13] Interest in the latter problem had originated in the 18th century, when the solvability by radicals of an $n^{th}$ degree equation had been connected to the number of values a resolvent function of $n$ variables could take. But even though such issues had originally been closely related to equations, they would give rise to autonomous developments on substitutions, permutations and arrangements. Augustin-Louis Cauchy (1815, 1844-46), Joseph Bertrand, (1845), and Serret (1849), had especially stated results on some boundaries to the number of values of functions. At the turn of the 1850s-1860s, substitutions themselves would become the main focus of the series of papers published by authors such as Émile Mathieu, Jordan, or Thomas Kirkman [Ehrhardt 2007, p. 291-393].

Jordan's approach consisted in reducing the problem from the class of transitive *n*-ary substitutions (corresponding to irreducible equations) to the class of primitive substitutions. As in Galois's 1830 works (that Jordan had not studied yet), this reduction was modeled on Gauss's method of successive factorizations of (irreducible) cyclotomic equations $x^{p-1}=1$, $p$ prime. The factorizations resulted from organizations of the roots in a specific order by appealing to the two indexings provided by a primitive root of unity and by a primitive root mod. $p$ [Neumann 2007]. The roots were then decomposed into "periods", or "groups" as Poinsot had designated them in 1808. This decomposition resorted to a single kind of substitutions (i.e. cycles). But two forms of actions had to be distinguished depending on whether the cycles were acting within the groups or between the groups. Poinsot had discussed these two forms of actions from a geometric perspective. The roots generated by a primitive root of unity could be represented "as if they were in a circle" [Boucard 2011, p. 68]. They could then be made to move forward by translations, i.e., by the operation (*i i+1)* on their indices. But they could also be made to move by rotations of the full circle i.e. (*i gi*). This analytic representation of substitutions played a crucial role in Galois's approach, as well as in Enrico Betti's and Jordan's ulterior works [Brechenmacher 2011a].

---

[13] This problem is tantamount to finding the possible orders for subgroups of the symmetric group. Given a function $\varphi(x_1, x_2,...,x_n)$ of $n$ "letters", a "value" of $\varphi$ was a function obtained by permuting the variables, *i.e.*, for any $\sigma \in Sym(n)$, $\varphi^\sigma(x_1,x_2,...,x_n)=\varphi(x_{1\sigma},x_{2\sigma},...,x_{n\sigma})$ was a value of $\varphi$. If $\varphi$ takes only one value, then it is symmetric and can therefore be expressed as a rational function of the elementary symmetric functions. If $x_1,...,x_n$ are the roots of an equation with coefficients on a given "rational domain", this means that $\varphi$ can be expressed as a rational function on the rational domain. In general, $\varphi$ can take up to $n!$ distinct values and to intermediary cases between 1 and $n!$, normal subgroups of the symmetric group can potentially be associated to $\varphi$ by considering the set of substitutions leaving $\varphi$ invariant. If $\varphi$ takes, for instance, $\rho$ distinct values $\varphi_1$, $\varphi_2$, ..., $\varphi_\rho$, these values can be considered as the roots of an equation of degree $\rho$ whose coefficients are the symmetric functions of the initial variables.

In 1860, Jordan followed Poinsot's reformulation of Gauss's decomposition in dividing the letters on which the substitutions are acting into "groups," each of a same cardinal $p^n$, $p$ prime.[14] The system $T$ of substitutions was then simultaneously partitioned into a "combination of displacements between the groups $\Gamma_i$ and of permutations of the letters within each of the groups" [Jordan, 1860, p. 5].[15] As a result, the substitutions of $T$ were divided into two "species". On the one hand, inside each block $\Gamma_i$, the letters were cyclically substituted by first specie of substitutions (*i gi*) operating on the powers of a primitive root mod.*p.* On the other hand, second specie of substitutions (*i i+1*) substituted cyclically the blocks themselves. Each specie of substitution corresponded to one of the two forms of representation of cycles. Their products generated linear forms (*i ai+b)*; $T$ was therefore what Jordan would designate in the mid 1860 as the *linear group*.[16]

In regard to the contemporary works of Charles Hermite, Joseph-Alfred Serret, or Mathieu, the originality of Jordan's approach was its focus on some general properties of classes of *n*-ary substitutions. This form of generality contrasted with the one of the explicit exhaustive lists of analytical forms of substitutions acting on 5 or 7 letters Hermite had given in connection to his works on the reduction of the degree of the modular equation of order 5, a problem that was also related to the legacy of Galois [Goldstein, 2011]. It also contrasted with the *general* results that authors such as Cauchy and Bertrand had stated for *all* substitutions.

Unlike these two forms of generalities, Jordan had appealed to the framework of the theory of order and to its focus on the relations between classes of objects. When the number of values of a function was less than *n!*, he considered that a "symmetry occurred within the function" as an application of "what Poinsot has distinguished from the rest of mathematics as the theory of order" [Jordan 1860, p. 3]. Even though he acknowledged his method of reduction from a general class of group to a less general, simpler, one, was not efficient for applications, Jordan claimed that in the aim of "studying the problem of the symmetry in itself, the method is not only more direct, it is also more natural and is actually the only way that leads to the true principles" [Jordan 1860, p. 4]. Jordan also highlighted the analogy between his method and the reduction of a helicoidal motion into motions of translation and rotation. Moreover, he eventually claimed that what one may designate as the unscrewing of groups was the "very essence" of the question [Jordan 1860, p.5].

Further developments of the "method of reduction" would give rise to several of Jordan's key results. After 1864, Jordan investigated the reduction of the general linear group into chains of normal subgroups. It was in this context that he stated the theorem on the invariance of the length and of the composition factors of the compositions series of a group, i.e., what is nowadays designated as the *Jordan-Hölder theorem*.[17] It was also in this context that Jordan introduced in 1868 what would later be designated as *Jordan's canonical form* of linear substitutions [Brechenmacher 2006a]. Later on, the method of reduction would eventually lie beneath the architecture of the 1870 *Traité* [Brechenmacher 2011a].

---

[14] The ambivalence of the terminology "group" as regard to the distinction between the "permutations of the roots" and the "substitutions" was the very nature of "groups" as they originated in Poinsot, Galois, or Jordan's investigations on the decomposition of imprimitive groups by the consideration of blocks of imprimitivity of letters.

[15] From a modern perspective, the "groups of permutations" correspond to a decomposition of the field into blocks of imprimitivity $\Gamma_i$ under the action of an imprimitive substitution group $T$, which is itself decomposed into a primitive quotient group. See [Neumann, 2006].

[16] In modern parlance, this group is actually an affine group.

### 1.3. The *irrationals (*1868-1870):
### modular equations, the 27 lines and elliptic/hyperelliptic functions

From 1864 to 1868, Jordan gradually transferred to the reduction of groups the "essence" he had originally attributed to the theory of order. The series of papers he published in 1868-1869 were fitting into a new thematic organization grounded on substitutions and equations.

Jordan claimed he was aiming at a "higher point of view on the transformation and the classification of the irrationals" [Jordan, 1870, p. V]. Jordan's "*irrationnelles*" pointed to both quantities and functions. They implicitly referred to the works of authors such as Hermite, Kronecker, Betti, and Brioschi who had considered that the impossibility to solve general algebraic equations of a degree higher than five raised the issue of identifying the nature of the "orders of irrationalities" associated to general equations of higher degree. This program involved finding the "most general" functions by which the roots of a general equations of a given degree could be expressed, such as in the solutions Hermite and Kronecker gave to the general quintic by using elliptic functions.

But in contrast with Jordan's, Hermite's and Kronecker's approaches to the "orders of irrationalities" were fitting into the research field of *arithmetic algebraic analysis* that had developed between the 1820s and 1850s [Goldstein and Schappacher 2007a, p. 26]. The minor role this research field attributed to geometric aspects – even the ones connected to Gauss's *Disquitiones arithmeticae* – highlights that it was different in nature in regard with Jordan's theory of order, i.e., with what has been designated above as the *regular solid* research field. As a matter of fact, the unity of the research field of arithmetic algebraic analysis had been torn apart in the 1860s. In contrast, geometric approaches to elliptic functions through invariant theory such as Clebsch's were playing a more and more important role. The evolution of Jordan's works in the 1860s gave rise to one of the lines of developments that participated to the dismemberments of transversal organizations of knowledge and to the emergence of object-oriented disciplines.

This evolution is exemplified by the "Commentaires sur Galois*"* that Jordan published in *Crelle's journal* in 1869. This paper laid the ground for the opening chapter of the third section of the *Traité*: "*Théorie générale des irrationnelles.*" Indeed, Jordan claimed that the interplays between Galois adjunctions of roots "to an equation" and the resulting reduction of a group could be "very useful for the classification of algebraic irrationals" because the invariance of the number of links in the so called Jordan-Hölder reduction of a group "gives a very distinct definition of the "degree of irrationality" of the roots of a given equation" [Jordan 1866, p. 1064].

---

[17] The procedure of chain reduction Jordan began to resort to would not only develop later into the Jordan-Hölder theorem but also into a false theorem stating successive reductions of blocks of imprimitivity [Neumann, 2006, p. 413].

But most of the third section of the *Traité* was actually devoted to applications of Galois theory. Moreover, while the "algebraic applications" mainly presented classical results on the solvability of cyclotomic, abelian, and metacyclic equations, it was the geometric and transcendental applications that constituted the actual core of the section "*Des irrationnelles.*" This part of Jordan's book presented an original synthesis of the legacies of the works of Alfred Clebsch and Hermite on the special equations associated to elliptic/hyperelliptic functions. Most of the transcendental applications were indeed devoted to the modular equations and to the equation of the division of the periods of elliptic functions. Jordan especially presented in the framework of substitutions groups Hermite's solution of the general quintic through the modular equation of degree 5 [Goldstein, 2011]. Jordan concluded this section by proving that the general equations of degrees higher than five cannot be solved by elliptic functions, These equations therefore call for investigations of hyperelliptic functions of multiple periods.

In this perspective, Jordan appealed to Clebsch's geometric approach to the invariants of binary forms. On the model of Galois's, Betti's and Hermite's reduction of the degree of the modular equation of 5, and of Clebsch and Paul Gordan's works on the bisection of the periods of abelian functions, Jordan investigated the reduction of the degree of the equation of the trisection of the periods of hyperelliptic functions of four periods. In 1869, he discovered that the reduction of this $80^{th}$ degree equation lead to a group he recognized as the one of the equation of the 27 lines on a cubic surface that had been discovered by Cayley and Samon, and that had been investigated by Steiner. It was at this occasion that Jordan announced the forthcoming publication of the *Traité.*

The theorem connecting the 27 lines and hyperelliptic functions played a crucial role in the *Traité*. It concluded the geometric applications of the section "*Des irrationnelles,*" which it connected to the applications to transcendental functions. Moreover, the 27 lines actually supported Jordan's claim that the substitutions groups introduced by the invariance of some algebraic forms permitted "an investigation of the hidden properties of the equations [of geometry]" [Jordan 1869b, p. 656] such as the (already known) fact that the group of the 28 double tangents (*Sp₆(2))* can be reduced successively to the groups of the 27 lines on a cubic (*Sp₄(3))* and of the 16 straight lines on a quartic surface having a double conic (*Sp₂(2)).*[18]

In sum, the *Traité's* third section claimed the relevance of substitutions groups for developing a higher point of view on the "order of irrationalities" related to the special equations of elliptic / abelian functions. It therefore inscribed substitutions in a framework completely different than the one of the theory of order to which Jordan had originally appealed. While no consideration on polyhedrons, crystallography, *analysis situs*, or groups of motions appeared in the *Traité*, Jordan attributed an important role to Clebsch's geometric approach and to Hermite's works on modular equations.

---

[18] The symplectic group $Sp_{2n}(p)$ was introduced by Jordan in 1869 as leaving invariant a non degenerate alternate bilinear form.

In the late 1870s, Jordan returned to issues related to special equations in investigating the hypoabelian group (i.e. an orthogonal group in a field of characteristic 2) attached to theta functions of four variables in connection to the general equation of the eight degree, thereby following some recent works by Max Nöther and Heinrich Weber.

### 1.4. The aftermath of the *Traité*: the limit of transitivity of substitution groups (1870-1874)

In addition to Jordan's punctual interest for probability in connection to his teaching at the *École polytechnique* (as a substitute for Hermite) and to a few papers connected to earlier concerns for motions of solid bodies or for the topology of lines and surfaces, most of Jordan's works in the early 1870s were devoted to investigations on classes of substitutions groups.

Jordan nevertheless did not refer to Galois any longer after 1870. It is therefore likely that the *Traité*'s emphasis on the issue of the solvability of equations was at least partly aiming at legitimating the general approach to classes of *n-ary* substitution groups that Jordan had developed since his thesis. In the early 1870s, Jordan indeed returned to issues related to the problem of the number of values of functions in the tradition of the boundaries Cauchy, Bertrand, and Serret had stated for such number of values. He especially appealed to the approach Mathieu had developed in the early 1860s on multiply-transitive groups. In 1875, Jordan eventually stated two finiteness theorems on the limit of transitivity of a substitution group and on the limit of the number of "classes of primitivity" of a group.[19]

In parallel, Jordan continued to advocate the essential nature of substitution groups in applying to new topics the methods he had developed for dealing with general linear groups - especially the canonical reduction of linear substitutions, e.g., linear differential equations with constant coefficients, the mechanical stability of small oscillations, and the theory of bilinear and quadratic forms. Jordan's general *n* variables approach to linear groups also laid the ground for a memoir on $n^{\text{th}}$ dimensional geometry (1875).

In the mid 1870s, Jordan's focus on relations and classes of objects was explicitly criticized by Kronecker for its false generality and formal nature [Brechenmacher 2007]. The latter accused Jordan of having mixed up tools relative to the orientation he had given to his investigations (i.e. *n-ary* linear substitutions) with the inherent significations of "objects of investigation." Moreover, he accused Jordan's *n* variables approach to resort to a false generality because of the non-effectiveness of its methods that appealed to solutions of general algebraic equations of degree *n*. After the 1874 controversial episode with Kronecker, Jordan stopped working on general classes of groups. From then on, he rather developed a group theoretical approach to specific object of investigations.

---

[19] Jordan calls the "class" of a group *G* the smallest number $c>1$ such that there exists a substitution of *G* which moves only *c* objects. The finiteness theorem states that there is an absolute constant *A* such that if *G* is primitive and does not contain the alternating group, then $n \leq Ac^2 \log c$ (in other words, there are only finitely many primitive groups of given class *c* other than the symmetric and alternating groups).

## 1.5. Differential equations and algebraic forms (1874-1881)

In the mid 1870s, Jordan returned to the notion of group of monodromy to which he had devoted his second thesis in 1860 in the legacies of Cauchy, Puiseux, and Hermite. He then appealed to some recent works by Fuchs, Hamburger, and Frobenius, to reduce the problem of the algebraic integration of Fuchsian linear differential equations of degree $n$ to the one of the classification of the finite linear subgroups of $Gl_n(\mathbb{C})$ [Brechenmacher 2011b]. From 1876 to 1880, Jordan provided such classifications for $n=2,3,4$.[20] Moreover, he proved than any periodic linear substitution can be reduced to a diagonal form. The latter result was then used as a lemma to prove Jordan's third finiteness theorem on the index of any finite linear group relatively to a normal abelian subgroup.[21] This group-theoretical approach to linear differential equations was completed by a series of papers that appealed to Gordan's invariants and covariants of binary forms.

Complementary to his investigations on linear differential equations, Jordan developed a specific approach to algebraic forms in Hermite's legacy. Jordan's series of papers on algebraic forms mixed some typical arithmetic-algebraic-analytic Hermitian concepts [Goldstein, 2007], with his own practice of canonical reduction of linear substitution. More precisely, the articulation of Jordan's works on finite linear groups with Hermite's notion of reduction of quadratic forms resulted in the statement of a fourth finiteness theorem for unimodular classes of equivalences of algebraic forms of degree $n$ of non vanishing determinant.[22] Hermite's continual reduction also gave rise to a method of decomposition of infinitesimal transformations in the framework of Lie's continuous group theory.

## 1.6. Collective responsibilities

The nature of Jordan's works changed in the early 1880s when the latter accessed to most of the key positions of Parisian mathematics: professor at the *École polytechnique* (1876), elected member of the geometry section of at the *Académie des sciences* (1881), professor at the *Collège de France* (1883), and director of the *Journal de mathématiques pures et appliquées* (1885).

From 1884 to 1904, Jordan published a number of academic notes, reports, and obituaries. In connection to the first edition of his *Cours d'analyse* [1882], he extended to functions of bounded variation the Dirichlet condition for the convergence of the Fourier series (1881). Later on, in the second edition of his *Cours* [1892], he gave a definition of simple and multiple integrals, i.e., of the measure of a curve. Jordan also proved his famous theorem on the condition under which a curve divides the plane between the points in its interior and the points in its exterior [Guggenheimer 1977].

---

[20] Jordan's 1876 classification of the finite subgroups of $Gl_2(C)$ was nevertheless incomplete as was pointed out by Klein. The 1878 classification of finite subgroups of $Gl_3(C)$ had also missed a class of simple group as was pointed out by Valentiner in 1889. See [Dieudonné 1962].

[21] For any integer $n$, there exists a function $\varphi$ $(n)$ such that any finite group $G$ of matrices of order $n$ contains a normal subgroup $H$ which is conjugate in GL $(n \ C)$ to a subgroup of diagonal matrices, and such that the index $(G: H)$ is at most $\varphi$ $(n)$ (equivalently the quotient group $G/H$ can only be one of a finite system of groups, to isomorphism).

[22] In nowadays parlance, Jordan considered the vector space of all homogeneous polynomials of degree $m$ in $n$ variables, with complex coefficients; the unimodular group $SL$ $(n, \mathbf{C})$ operates in this space, and Jordan considered an orbit for this action. Within that orbit he considered the forms having (complex) integral coefficients and he placed in the same equivalence class all such forms which are equivalent for unimodular substitutions having (complex) integral coefficients. The number of these classes is finite, provided that $m>2$ and that the discriminant of $F$ is not zero.

During this period, Jordan published a few papers in response to some other authors' interest for his 1870s investigations in group theory, such as the works of Bochert and Maillet on the limit of transitivity / the degree of primitivity of a group (1895, 1908), the works of Frobenius and Hölder on solvable groups (1898), or the ones of Dickson on linear groups in Galois fields (1904). After 1904, Jordan eventually returned to investigations of forms (or networks of forms) invariant by a given substitutions group. He also returned to some issues he had investigated at the beginning of his career in connection with the *regular solids* research field, as for the topological considerations on constellations of points (that may represent a polyhedron) that laid the ground for Jordan's last paper.[23]

## 1.7. Conclusions

In sum, Jordan's works passed through several collective organizations of knowledge from the 1860s to the 1920s. The works of the 1860s highlight the legacies of some transversal fields of research such as the ones of the *regular solids* and of *arithmetic algebraic analysis*. But these works also point to some more local traditions articulated to these two European-scale research fields, such as Cauchy's approach to the problem of the number of values of functions or Poinsot's theory of order. In the mid 1860s Jordan had aimed at a new unified approach to the "essence" of the "irrationals" by articulating these legacies to some more recent works, such as Hermite's and Clebsch's. The approach to substitution groups he developed in his *Traité* eventually presented an original synthesis between these influences. Here, Jordan's reference to Galois was crucial for bringing together various issues such as the number of values of functions, the solvability by radicals, the special equations of geometry and of elliptic/hyperelliptic functions, etc., and for binding them all by articulating substitutions to equations.

In contrast with organizations such as the theory of order, the *Traité*'s synthesis thus emphasized and object-oriented organization of knowledge. On the one hand, the theory of order was torn apart as some of its key concepts, such as groups or cyclotomy were separated from others such as symmetry, polyhedron, motion etc. But on the other hand, the essential nature Jordan attributed to his "method of reduction" of the analytic representation of linear groups had underlying it some specific articulations grounded on cyclotomy in the legacy of the theory of order. The process of reduction of a group into chains of subgroups was for instance implicitly appealing to the unscrewing of a helicoïdal motion into some rotation and translation motions on the model of Poinsot's interpretations of cycles as two forms of motions on a circle.[24]

---

[23] This topic points to the tradition of the topological approach to the deformation of networks of points in connection to crystallography [Scholz 1989, p. 93].
[24] On the connections between Poinsot's 1851 investigations and crystallography [Scholz 1989, p. 93].

Because Jordan's "method of reduction" was actually lying beneath the architecture of the *Traité*, this book presented both an object-oriented theory and a transversal organization of knowledge. Various forms of readings of Jordan's book would therefore develop.[25]

In the 1870s, the attribution of an essential nature to relations between classes of objects was laying the ground for Jordan's approach to both the contemporary mainstream issue of Fuchsian linear differential equations and to the more local legacy of Hermite's theory of forms. Later on in the 1880s and 1890s, Jordan's publications were mainly inscribed in the discipline of Analysis as it was taught at *École polytechnique*, and which was usually considered in France as grounding the unity of mathematics [Brechenmacher 2012].

## 2. Individual and collective forms of references to Jordan (1868-1894)

Let us now investigate the corpus of the *Jahrbuch* reviews that referred to Jordan from 1868 to 1940. The diagrams n°2 give an overview of the distribution of the Jordan corpus. The evolutions shown by the first diagram have to be considered in the perspective of the parallel exponential growth of mathematical publications. The corpus related to the five-year time-period 1931-1935 is as big as the one related to the forty-year period 1870-1910. The proportions represented in the second diagram show a relative decline after 1884 when Jordan's publication production decreased. At the turn of the century the relative weight of the Jordan corpus in the *Jahrbuch* increases abruptly. From 1910 to 1930, the Jordan corpus keeps increasing at a faster pace than the global production of mathematical papers. In the 1930s, the Jordan corpus eventually represents about 5-6 $^0/_{00}$ of the whole *Jahrbuch*.

### 2.1. An overview of the periodization

A finer analysis of the global corpus results in a partition into two groups of three time-periods: ([1868; 1877]; [1878; 1884]; [1885; 1894]); ([1895; 1909]; [1910; 1929]; [1930; 1939]). As we shall see, specific references to Jordan lie beneath each of these time-periods. The evolutions of the first group of time-periods mainly resulted from the impulses Jordan himself gave to his own works. In contrast, most post-1895 evolutions have underlying them some collective trends mainly independent from Jordan.

The time-periods of the first group are close to the ones considered in the previous section when commenting on Jordan's works (recall that there is no *Jahrbuch* database for the period prior to 1868). Actually, the time-periods considered in this paragraph mainly differ from the ones considered in the previous section by three-year intervals. These intervals may be understood as times during which Jordan's papers were commented on as newly published works, i.e., as times of contemporaneity to Jordan's publications.

The first time-period is dominated by references to the *Traité*. It therefore mixes allusions to the *regular solids* research field, to the traditional problem of the number of values of functions, and to the equations associated to elliptic and hyperelliptic functions. The distribution of the *Jahrbuch* mathematical classifications in the Jordan corpus (diagrams n°3) is close to the one of Jordan's own works (diagrams n° 1) except for the chapter "equations."

---

[25] On the historical notion of readers of a text, see [Goldstein 1995].

In contrast with the fact that Jordan's works were often classified in the latter chapter, only a very few reviews referred to Jordan in the chapter "equations," even on a larger time-scale than 1868-1877. As shall be seen later, the reception of the *Traité*'s presentation of Galois theory sheds some light on this situation.

The period 1878-1885 is dominated by the issue of the algebraic integration of linear differential equations: it starts with Fuchs's reference to Jordan's "method of substitutions" and ends when the latter method was diluted in the approaches developed by actors such as Klein, Poincaré, or Picard.

The period from 1885 to 1894 witnesses an increasing variety in the distribution of the *Jahrbuch* classifications. This diversification is mainly due to the subsistence of earlier uses of Jordan's works in addition to the appearance of new forms of references, mainly to Jordan's *Cours d'analyse* .

The second group of time-periods comes with a quantitative shift in the Jordan corpus. This evolution is mainly due to the references that were made to the *Traité*'s linear groups in Galois fields by a flood of texts published by Dickson and his followers in Chicago.

The explosive growth of references to Jordan in linear group theory was nevertheless short lived. The growth of the Jordan corpus after 1910 is indeed mainly due to the universalization of some generic designations of concepts/theorems/methods, i.e., the Jordan curve theorem, the Jordan measure, and the Jordan-Hölder theorem after 1910, as well as the Jordan canonical form theorems after 1930.

In looking more closely at each of the six time-periods, we shall especially aim at describing some collective types of references to Jordan in their chronological order of appearance.

## 2.2. Forms of distributions of references in time

Before analyzing further the Jordan corpus, it is compulsory to discuss first the various forms of references to Jordan that can be found in this corpus. Several different aspects could be analyzed in regard with the ways texts are referring one another, e.g., whether the reference is a citation, a comment, a refutation, a simplification, a generalization, etc. But such aspects are not much relevant for the overview on the whole Jordan corpus we give here. The present paper indeed rather aims at investigating some forms of distributions of references in time. We shall distinguish between four forms of such distributions: punctual references, continuous references, sporadic references, and generic references.

The *punctual form* corresponds to a specific reference shared by a very few papers that were being published during a limited time-period. For instance, the references to Jordan's 1875 approach to probabilities have a punctual form within the Jordan corpus (Laquière 1880; van der Berg, 1891).

The *continuous form* corresponds to references that are explicitly appealing to a broader collective organization of knowledge, such as a discipline or a research field. For instance, Jordan's works on invariants and covariants were discussed in connection to the traditional British approach to invariant theory, mostly in Great Britain (Sylvester 1881; Forsyth 1882;

Grace 1903; Young, 1903; Wood, 1904) but also in the United States (Mac Kinnon 1895; Keyser 1898).

The *sporadic form* contrasts with both forms of references introduced above. It points to series of punctual papers that refer one to another, thereby giving rise to a network of texts with a specific identity. In some cases, Jordan's works themselves can carry on this specific identity, such as in the case of the references to Jordan's canonical from 1870 to 1900. Some forms of references to Jordan were sporadic right from the beginning, such as the reference to the 1875 paper on the geometry of $n$-dimensional Euclidean spaces (Lipschitz 1875; d'Ovidio 1877, Brunel 1881; Cassani 1885; Biermann 1887; Predella 1889; Lovett 1901, 1902). As shall be seen in greater details later, the fate of most *continuous forms* of references to Jordan's works was to turn either sporadic or generic at some point, even when the collective organization that was originally involved in the reference used to be a very large one.

The *generic form* is one the possible fate for both *continuous* and *sporadic* references. A reference is generic when it does not point to any specific part of Jordan's work (even implicitly). The point here is not so much about whether or not it is possible to find out which specific parts of Jordan's work is implicitly connected to a generic reference. The point is that the very nature of a generic reference is to point (at least implicitly) to several works at the same time without specifying any of these. For instance, references to the Jordan curve theorem turned generic at the beginning of the 20[th] century. Such references then implicitly pointed to some works of Schönflies or Osgood as well as to the ones of Jordan. We shall see later that the references to Jordan canonical form theorem give an example of a sporadic form that first turned generic in some local contexts at the turn of the century before it eventually became generic at a global scale in the 1930s.

### 2.3. Early references to Jordan (1868-1877)

The distribution of the languages in the 26 reviews that were mentioning Jordan from 1868 to 1877 (diagrams n°4) points to the role played by Jordan's own network of correspondence. In the late 1860s, Jordan was indeed mostly in epistolary contact with Italian (Brioschi, Cremona) and German (Clebsch, Borchardt) authors on issues related to the equation of the 27 lines on a cubic.

The domination of the German language is nevertheless not a trivial situation in the context of the aftermath of the 1870 Franco-Prussian war. A few testimonies actually presented as a matter of national pride both Jordan's claim to follow the legacy of Hermite's works on the general quintic and Jordan's capacity to deal with some contemporary German researches. In 1866, Bertrand had already highlighted the importance of Jordan's works in connection to his celebration of Hermite's achievements in his national *Rapport sur les progrés les plus récents de l'analyse mathématique.* Moreover, in 1873, Max Marie alluded to the loss of Alsace-Lorraine when commenting on some of Jordan's results on equations solvable by elliptic functions: "this, at least, has been taken back from the Germans" [Marie 1873, p. 943].[26]

### 2.3.1. Substitutions

Jordan's *Traité* has often been described as the first autonomous presentation of both group theory and Galois theory. As already alluded to before, the reception of the book has often been considered to have been very limited, either because of the novelty of its group theoretical approach, or, on the opposite, because of its outdated focus on substitutions groups at a time when Klein and Lie were developing researches on transformation groups [Wussing 1984].

Jordan's theory of substitutions was nevertheless immediately presented to the audience of *Battaglini's journal* by Janni, one of the editor's students [Martini 1999]. But both Janni's series of four papers, and the presentation Sardy gave in the same journal of Jordan's theory of higher congruences show that the novelty of Jordan's approach to substitutions groups was more related to its synthetic nature than to the presentation of a genuine new theory.

That the *Traité* was considered in the continuity of previous works is made clear by the series the papers in which Netto reformulated Jordan's approach in the legacies of Cauchy and Kronecker. Like Janni, Sardi, or Netto, most auhors who referred to the *Traité*'s general presentation of substitutions groups pointed to the synthetic nature of the book, which they usually discussed in regard with the works of authors such as Cauchy, Kronecker, Hermite, Bertrand, Serret, and Mathieu. I shall designate this form of reference to Jordan as the *synthetic-Traité* type. It can be found in ulterior time-periods such as in (Netto 1874, 1877, 1878, 1882), (Pellet 1887), (Bolza 1890), (Borel et Drach 1895), (Vogt 1895), (Weber 1895), (Picard 1897), (Echegaray 1897), or (Pierpont 1900).

In connection to the *synthetic-Traité* type of reference to Jordan, the latter's specificity was often attributed to the generality of his approach to linear groups, which most authors actually rejected. In 1882, Netto grounded his own treatise on substitutions on an ideal of effectiveness in the legacy of Kronecker as opposed to the "abstraction" of Jordan's developments on *n*-ary linear groups [Brechenmacher 2007]. As shall be seen below, some specific aspects of the *Traité* had nevertheless circulated early on

Recall that for decades, both Galois's legacy and the emergence of linear group theory opposed two approaches which both aimed at reaching the "essence" of mathematics [Brechenmacher 2011a]. On the one hand, some authors, following Hermite and Kronecker, aimed at characterizing the *special* nature of *general* equations of a given degree. On the other hand, some other authors, following Betti and Jordan, focused on the *relations* between *classes* of solvable equations (or groups) of a *general* degree *n*.

The two approaches were nevertheless both presented in Jordan's *Traité.* The *Traité*'s third section, "*Des irrationnelles*," presented Jordan's own synthesis on the first approach to general equations. As shall be seen in greater details later, from 1870 to 1890, Jordan's synthesis would mainly be considered in the continuity of previous works and its specificity would again be rejected. For this reason, until the turn of the century, "linear groups" usually designated the groups of binary or ternary unimodular fractional linear substitutions (i.e. $PSL_2(p)$ and $PSL_3(p)$) Kronecker, Klein, and their followers had investigated in the legacy of Galois's, Betti's, and Hermite's approaches on "the three Galois groups" of the modular equations of order 5, 7, and 11 [Goldstein 2011]. The second approach to general equations went with Jordan's specific method of reduction. It structured the *Traité* in a complex chain of generalizations of special model cases. As shall be seen later, Jordan's focus on general *n*-ary linear groups had almost no circulation until this approach would be developed in a network of text revolving around Dickson's works.

---

[26] On the echoes of the war in the French mathematical community in the 1870s, see [Gispert 1991].

### 2.3.2. Special equations

Unlike the *synthetic-Traité* type of reference to Jordan, some *Jahrbuch* reviews pointed to specific innovations of Jordan's *Traité*. These did not refer to the general presentation of substitutions groups but to the *special equations* related to elliptic/hyperelliptic functions. In his review of the *Traité*, Netto especially highlighted the theorem in which Jordan had discussed the possibility to solve general equations of degree higher than five by hyperelliptic functions. This theorem was actually the only result Netto's review explicitly quoted. Moreover, it was the theorem that connected the 27 lines on a cubic to hyperelliptic functions that Cremona cited when he addressed his congratulations to Jordan for the publication of the *Traité*.

As was said before, issues related to elliptic/hyperelliptic functions had traditionally been considered in the transversal research field of arithmetic-algebraic- analysis from the 1820s to the 1850s. In the 1870s, papers related to the special equations of elliptic/hyperelliptic functions were still usually classified in various sections, especially "2.1. Equations," "2.2. Theory of forms" (Clebsch 1869), "2.3. Substitutions," "4.5. Analytic geometry" (Geiser 1869), "7.1. Function theory, generalities" (Marie 1873), and "7.2. Special functions." Such classifications reflected before all the choices of the *Jahrbuch* reviewers. But these choices had nevertheless underlying them some patterns of organizations of knowledge. To be sure, Jordan's repeated claims to aim at focusing on equations/substitutions made it likely for his papers to be classified accordingly. This situation sheds light on the contrast between the role played by the section "equations" in the Jordan corpus in regard with its role in the classifications of Jordan's own paper.

But this situation also highlights the coexistence of some traditional transversal organizations of knowledge with several newer object-oriented theories. As a matter of facts, most papers that refered to Jordan's special equations in the 1870s aimed at reformulating the *Traité*'s focus on substitutions groups. Jordan's results on the relations between the 16 lines on a quartic, the 27 lines on a cubic, and the trisection of the periods of abelian functions were reformulated in the framework of Clebsch's geometric approach to the invariant theory of binary forms (Geiser 1869; Clebsch 1869; Brioschi 1877),[27] while Jordan's approach to the division equation of elliptic functions was reformulated by Sylow (1871), Marie (1873), and Kronecker (1877). Most authors highlighted the specificity of Jordan's "method of substitutions," even though they usually avoided this method by appealing to other approaches to elliptic functions, such as the ones of Weierstrass and Leo Königsberger (Krause 1881 ; Nöther 1879, 1880, 1881).

Later on in the 1880s, Klein and his followers would interlace the *special equations* type of reference to Jordan with the one to the *regular solids* research field (Klein 1877; Gierster 1881; Hess 1883; Hurwitz 1884; Burkhardt 1893).

### 2.3.3. Regular solids

From 1868 to 1877, about 67% of the reviews of the Jordan corpus were pointing to the *Traité*. Most of the remaining third actually pointed to the 1868 memoir on "groups of motions." The papers under review here were mostly fitting into the *regular solids* research field in the legacies of authors such as Poinsot, Bravais, Bertrand, or Listing.

As for the case of the *special equations* type of reference, the *regular solids* type was transversal to the *Jahrbuch* classification. Here, the chapters involved were "8. Pure geometry. 2. The notion of continuity (analysis situs)" (Bertini 1869, on the symetries of polyhedrons ; Becker 1869, 1873, on the deformation of surfaces), "9. Analytic geometry. 3. Surfaces" (Boussinesq 1872), "10.Mechanics. 2. Kinematics" (Linguine 1872, on the motions of solids), "11. Mathematical physics. 1. Molecular physics" (Sohncke 1875, on crystallography),[28] and "7.2. Special functions" (Schwarz 1872, on Riemann surfaces, linear differential equations, and Gauss's hypergeometric series).

The unity of the *regular solids* research field was nevertheless already in the process of being torn apart. Some authors especially separated topological aspects, such as laces of integrations, trees of ramifications and Riemann surfaces, from other issues (Polignac 1880, Cayley 1882). More importantly, and as was already said before, Jordan's polyhedrons and groups of motions were mixed with the special equations of elliptic functions and with linear differential equations in the framework of Klein's transversal approach to the icosahedron [Klein 1884].[29] This context favored the circulation of some specific aspects of the *Traité*, such as the notion of isomorphism that Jordan had transferred from crystallography to substitutions groups (Klein 1878; Capelli 1878). In the 1890s, Klein would retrospectively emphasize the connections between his 1872 Erlangen program and Jordan's "approach to discontinuous groups" in connection to crystallography [Klein 1893].

Despite these new lines of developments, the research field of *regular solids* would nevertheless not disappear until the turn of the century. The problem of the determination of the surfaces that have the same symmetries as a given polyhedron was set in 1887 as the *Grand prix* of mathematics at the Paris Academy. At this occasion, Kirkmann made a claim of priority over Jordan's 1866-1868 works on polyhedrons, which he argued he had just discovered.[30] That this claim was published as an addendum to the solution Kirkman gave to a question in the *Educational times* highlights the public dimension of the problem set for the *Grand prix*. In his review of Kirkmann's note, Lampe ironically pointed out the self-contradiction in accusing another author of deliberate plagiarism while admitting such a little knowledge of the mathematical literature.

---

[27] A similar Clebsch oriented reference to Jordan was also developed by Klein and Lie conjoined works on the 16 double points on Kummer surfaces (their three papers nevertheless do not appear in the Jordan corpus). See [Rowe 1989].
[28] On the connections between Jordan's works on groups of motions and Sohncke's ones, see [Scholz 1989; p. 110-114].
[29] Note that this was the sole paper Klein ever published in Jordan's Journal.

This episode exemplifies the subsistence of the unity of the specific mix of concerns coming along with references to regular solids. From 1880 to 1894, several works were still referring to Jordan's groups of motions in connection to the long term legacies of authors such as Poinsot and Bravais (Godt 1881; Goursat 1887; Kirkmann 1887; Fedorow 1891; Hagen 1894). During this time-period, the *regular solids* type of reference, taken altogether with its subsequent evolutions with special equations and topological curves represents about one third of the Jordan corpus.

### 2.4. Invariance of algebraic forms on the action of a given finite group

The new references to Jordan's works that appeared from 1878 to 1884 were mostly associated to the increasing roles played by the chapters "6.5. Differential equations," "2.2. Theory of forms," and "3. Arithmetic 2. Theory of forms." Later on, from 1885 to 1894, the new references to Jordan were mainly related to the *Cours d'analyse*: functions of small variations (Kronecker 1885; Hölder 1885; de la Vallée-Poussin 1891), the definition of simple or multiple integrals (Mansion 1888; Pochhammer 1890; Schellenberg 1892), as well as some other topics mostly related to some issues of convergence of series (Starkow 1885; Gerbaldi 1891; Cayley 1892).

Let us now consider in further details how the issues related to the algebraic forms invariant on the action of a given finite group were connected to Jordan's works on linear differential equations of Fuchsian types, which might have been one the most influential contribution of Jordan until the turn of the century.

As has been seen above, some issues on classifications of finite linear groups lay at the root of Jordan's approach to differential equations. But Jordan's 1876 classification of the finite subgroups of $Gl_2(\mathbb{C})$ was incomplete as was immediately pointed out by Klein who claimed his priority over Jordan. In 1875, Klein had indeed classified binary linear fractional transformations by considering geometrically the binary forms they are leaving invariant. He had therefore appealed to Jordan's 1868 investigations on the groups of motions of polyhedrons, which he had considered in both the frameworks of Clebsch's geometric invariant theory and of Hermite's approach to the general quintic. In 1876, Jordan's new approach to differential equation raised Klein's interest for this topic and thereby instigated new connections between Schwarz's 1872 approach on the hypergeometic equation, Riemann surfaces, the general quintic (Klein 1877), invariant/covariants (Gordan 1877), modular functions (Dedekind 1877), and groups of monodromy of linear differential equations (Fuchs 1878). These new connections eventually gave rise to Klein's icosahedron theory, i.e., the theory of automorphic functions [Gray, 2000].

But we have seen that from 1878 to 1881 Jordan had generalized his finiteness theorem on finite linear groups to algebraic forms. It was actually through the prism of Hermite's legacy that the specific approaches Jordan had developed for dealing with linear substitutions circulated to Poincaré (1881, 1884). In addition to playing an important role in Poincaré's group-theoretic approach to Fuchsian functions, Jordan's "method of reduction" deeply influenced Poincaré's algebraic practices. The latter then played a direct role in the circulation of Jordan's works in France by orienting Léon Autonne's doctoral works on Jordan's approach to the algebraic integration of Fuchsian equations (Autonne 1884) [Brechenmacher, 2011b].

As a result, the circulation of Jordan's "method of substitutions" in France mainly went with the legacy of Hermite's theory of forms in the framework of linear differential equations. This circulation nevertheless went underground: apart from some sporadic references of authors such as Poincaré, Autonne, Picard, or Goursat, most papers that resorted to the methods Jordan had developped made not explicit reference to the latter. For instance, even though Jordan had published in 1874 a memoir in which he had grounded the algebraic integration of Fuchsian equations on group theory, Jordan's works was neither referred to by Lie's 1873-1875 programmatic papers on a theory of continuous groups of transformations (or by later developments by this author) nor by Picard's differential Galois theory [Archibald, 2011]. Later on, Jordan's works were nevertheless considered as belonging to the framework of differential Galois theory (Epsteen, 1902).

In contrast with the lack of explicit recognition for Jordan's *general* approach to linear substitutions, some *specific* results were often referred to altogether as "Jordan's methods" in the context of the algebraic integration of differential equations, e.g., Jordan's classification of finite linear groups, the various finiteness theorems, and the characterization of periodic linear substitutions. The investigation of algebraic forms invariant by some given linear groups gave rise to the first generic form of reference to Jordan (Goursat 1885; Bolza 1893). Investigations of the like were indeed designated as "Jordan's investigations on the group concept" in Meyer's review of [Minkowski 1885]. Moreover, Netto alluded to a "method of proof of a Jordan's kind" in reviewing [Picard 1887,] while Hamburger and Wassilief referred to "Jordan's method of finite linear groups" when reviewing (Floridia 1884) and (Sawitch 1892) respectively.

Moreover, Jordan's investigations on the theory of forms circulated in connection to Hermite's legacy in number theory.[31] Issues on the invariance of algebraic forms for the action of a given group were indeed related to the long-term legacy of what Hermite had designated as the problem of "similar substitutions" after Lagrange's concept of "similar functions." [32] In the 1880s, Hurwitz inscribed Jordan's and Poincaré's works on algebraic

---

[31] On Hermite's legacy, see [Goldstein 1999 & 2007], [Goldstein and Schappacher 2007b].
[32] Lagrange had developed the notion of "similar" functions in the 1770s (Cf. [van der Waerden 1985, p. 81]). Two functions *f* and *g* of the roots of a given equation are called *similar*, if all substitutions leaving *f* invariant also leave *g* invariant. It then

forms in the context of Klein's icosahedron while Minkowski (1885, 1887) developed the geometric interpretations of the reduction of quadratic forms in networks of the plane, as investigated by Poincaré in the legacy of Hermite in 1870s as well as by Jordan in the 1860s in connection to crystallography. Other authors reactivated Jordan's works of the late 1860s in investigating Cremona's substitutions of contacts in connection to algebraic forms (Kantor 1885; Autonne 1886).

Bringing altogether differential equations and algebraic forms, the *algebraic forms* invariant by finite groups type of reference to Jordan eventually represented about 22% of the Jordan corpus in the time period 1885-1894.

### 3. The institutionalization of group theory (1895-1910)

As was said before, the important quantitative growth of the corpus from 1895 to 1910 is mostly due to the increasing number of texts on "group theory" that were published in English by American authors close to the "Chicago research school".[33] As a result, the continuous forms of references to *special equations* or to *regular solids* faded away even though Jordan's works remained an important reference in the context of the algebraic integration of linear differential equations (Fuchs 1896; Boulanger 1898; Epsteen 1902; Schlesinger 1906).

### 3.1. On the variety of references to Jordan's "groups"

The quantitative domination of the Chicagoan's references to Jordan must be handled with care. Practices of publications were indeed varied at the turn of the century. For instance, from 1895 and 1909, Dickson was the author of about half of the papers in English reviewed in the Jordan corpus, a number which amounted to the one of the total of papers in French or in German during this time period.

But various types of references to Jordan's works on substitutions nevertheless coexisted with Dickson's flood of papers on linear groups in Galois fields. In the context of the institutionalization of group theory, different types of references to Jordan actually gave rise to specific approaches to groups.

First, even though Jordan was rarely referred to in connection to *regular solids* anymore (Cotton 1900), crystallography was considered as one of the main roots of group theory as Schönflies made it clear at the Chicago congress of 1893. Second, and similarly, while direct references to special equations turned sporadic (Geiser 1896; Pech 1903), one of the main trends in group theory was animated by some followers of Klein whose concerns for collineation groups in projective geometry were in the continuation of older concerns for special equations (Valentiner 1889; Maschke 1895; Wiman 1896; Newson 1901; Bagnera 1905; Blichfeldt 1905, 1907). In this context, references to the 1877-1879 works on finite linear groups gave rise to specific approaches to the invariants of finite groups (Maschke 1895; Moore 1899) or to periodic matrices (Ranum 1907).

---

follows that *g* is a rational function of *f* and of the coefficients of the initial equation.
[33] See [Parshall 2004].

Third, another important sporadic form of references to Jordan pointed to the two finiteness theorems on limits of transitivity / degrees of primitivity of substitution groups. These two theorems had first given rise to continuous form of references in connection with the traditional problem of the number of values of functions (Netto 1877; Battaglini, Beltrami, Casoratti 1884; Rudio 1887) until they had been generalized by (Bocher 1887, 1889). They had then turned sporadic in the early 1890s, thereby giving rise to a specific trend in group theory (Maillet 1895; Miller 1894, 1897; Wilkonson 1899; Manning 1906-1914, Bieberbach 1911).

### 3.2. Linear groups in Galois fields

Let us now consider the main type of references to Jordan's works at the turn of the century, i.e., the references to linear groups in Galois fields. New interests for Jordan's *n*-ary general linear groups had risen in connection with Dickson's thesis on "The Analytic Representation of Substitutions on a Power of a Prime Number of Letters with a Discussion of the Linear Group."

As was said before, the analytic approach to *n*-ary substitutions through a "method of reduction" had been one of the main specificities of Jordan's works in the 1860s. But even though this specific approach gave its structure to the *Traité*, it would have a very limited circulation until the 1890s. Most group-oriented approaches, such as Klein's iscosahedron or Poincaré's theory of Fuchsian functions, were indeed dealing with binary (or ternary) linear fractional substitutions, which they considered as geometric transformations.

In order to represent analytically substitutions on a finite number *m* of letters as functions on integers mod. *m*, it is compulsory to index these letters. If *p* is a prime number, such an indexing is given by the roots of the cyclotomic congruence $x^p - 1 \equiv 0$ (mod. *p*). As had been shown by Galois, if *m* is taken to be the power of a prime $p^n$, the indices can be considered as the "imaginary numbers" solutions of the congruence $x^{p^n} - 1 \equiv 0$ (mod. *p*). Galois's number-theoretic imaginaries therefore made it possible to give an analytic representation of linear substitutions on $p^n$ variables. They had been used as such by Galois, Betti, Serret, Mathieu, and Jordan [Brechenmacher 2011a].

It was for a similar reason that Moore appealed in 1893 to the field $F_{p^n}$ - which he designated as the Galois field - in the aim of generalizing to substitutions on $p^n$ letters Klein's investigations on the group of unimodular binary linear fractional substitutions mod. *p*, i.e., for introducing the new system of simple groups *PSl₂(pⁿ)* as a generalization of the simple groups *PSL₂(p)*. But the novelty of Moore's 1893 investigations on *PSl₂(pⁿ)* had been questioned in regard with some of Mathieu's works of the early 1860s that Moore did not know until 1893. It was in this context that, between 1883 and 1896, Moore and his first doctoral student Dickson appealed to Jordan's *Traité* to investigate some of the works published in the 1860s by Mathieu, Hermite, Serret, and Jordan himself. All these works involved the analytic representations of substitutions on $F_p$ or $F_{p^n}$. The Chicagoan especially focused on Jordan's definition of the general linear group *Gl($F_{p^n}$)*, i.e. "Jordan's linear group," as Moore said, as the normalizer of an elementary abelian group (i.e. of the additive group attached to the ground field $F_{p^n}*$).

Let us pause to recall that Moore had first presented his paper on Galois fields at the occasion of the closing lecture at the 1893 Chicago mathematical congress. This paper has often been commented as one of the symbols of the influence of the abstract "characteristic of trendsetting German mathematics" on the emergence of both the Chicago research school [Parshall 2004] and the American mathematical research community [Parshall and Rowe 1994]. And, indeed, Moore obviously aimed at celebrating Klein who was the guest of honor of the congress. But recall that the Chicago congress took place during the World Columbian exhibition, which was the occasion of much display of national greatness, such as the German university exhibit and the series of lectures given by Klein as an official representative of the German government. [34] The fact that Moore's paper eventually resulted in the circulation of a specificity of Jordan's works that was foreign to Klein's investigations [Brechenmacher 2011a] highlights the difficult problems raised by identifying the scales at which some national, institutional, or local categories played or did not play a relevant role.

In the early 1890s, Moore and Dickson were not the only authors to appeal to Jordan's linear groups. Definitions of Galois fields similar to Moore's were given in [Borel &Drach 1895] as well as in (Burnside 1894). Moreover, several authors claimed independently to have abstractly identified $Gl(F_{p^n})$ as the group of automorphisms of an elementary abelian group (i.e., Frobenius, Hölder, Moore, Burnside, Le Vavasseur, and Miller). These parallel interpretations of Jordan's definition of the general linear group highlight that Jordan's *Traité* was still referred to for its synthetic nature at the time of the institutionalization of group theory.

In this context, Dickson eventually became the leader of a specific approach to general linear groups, which was devoting an important role to analytic representations and matrix computations. After he completed his thesis, Dickson indeed published floods of papers generalizing to $Gl(F_{p^n})$ various results Jordan had stated for $Gl(F_p)$ (these papers amounted to more than 20% of the Jordan corpus from 1895 to 1910).

Dickson's 1901 monograph on *Linear groups in Galois fields* was therefore often presented as a direct continuation of Jordan's *Traité*. Its focus on the analytic representation of *n-ary* linear group was indeed giving their unity back to issues that had been commented in various distinct frameworks in the 1870s and 1880s. Dickson especially followed Jordan's investigations on the special equations of elliptic/hyperelliptic functions through the one of some special subgroups of the general linear group, such as the hypoabelian group attached to theta functions (Dickson 1898) or the "new abstract simple group" of order 25920 related to the 27 lines and the trissection of abelian functions of four periods (Dickson, 1899, 1900, 1904).

---

[34] The great wheel of the exhibition was explicitly presented as aiming at challenging the Eiffel tower. Reflecting on national greatness, one may compare the failure of the German university exhibit [Parshall and Rowe 1994] with the success of the French bakery which location was printed on all entrance tickets to the fair.

### 3.3. The Dickson network

At the turn of the century, generalizing Jordan's results to Galois fields turned into a collective trend in the local context of the Chicago school. But the *Traité*'s method of reduction of the analytic representation linear groups was also laying the ground for a larger-scale network of texts. This collective of texts can nevertheless not be identified to a theory, such as group theory (Frobenius, for instance, did not take part in it), nor to a research school, such as the Chicago school (Maschke was not a part of it and, in a sense, Moore never adopted entirely Dickson's approach), or to a national framework such as the development of the American mathematical community.

It is therefore compulsory to investigate further the collective dimension of the references to Jordan's linear groups. From 1895 to 1909, the 60 papers in the Jordan corpus that are related to linear groups in Galois fields amount to about 50% of the corpus and to more than 90% of the texts classified in "group theory." I shall designate the network of intertextual references originating from this group of *papers* as the 'Dickson network'. This network is mainly constituted of texts that have been published by French and American authors between 1893 and 1907. This group initially involved actors in Chicago (Moore, Dickson, Ida May Schottenfels, Joseph H. Wedderburn, [35] William Bussey, Robert Börger, Arthur Ranum) and in Paris (Jordan, Émile Borel and Jules Drach, Raymond Le Vavasseur, Jean-Armand de Séguier, Léon Autonne) but quickly extended to actors in Stanford (George A. Miller, William A. Manning, Hans Blichfeldt), and to other individuals such as William.L. Putnam, Edward V. Huntington, or Lewis Neikirk.[36]

Dickson's works were at the center of most intertextual relations. It was through the light the latter thesis had shed on the relations between the works of Hermite, Jordan, and Mathieu that the notion of Galois field first circulated outside of Chicago (Miller 1897; Séguier 1901; Manning 1903; Jordan 1904).

In addition to Dickson's works, the collection of texts revolves around some shared references to some works of the 1860s: Hermite's works on the analytic representation of substitutions, Serret's works on Galois fields, Mathieu's works on multiply-transitive groups, and Jordan's works on linear groups. The Dickson network thus revolved around a two-fold periodization which was articulated by Jordan's *Traité.* It can actually be understood as the space of circulation of a specific relation Jordan had established to the works of Galois in the 1860s. This circulation especially concerned specific practices for dealing with the analytic representation of linear substitutions such as Jordan's canonical form (Dickson 1898-1908; Moore 1898; Burnside 1898; Bromwich 1899, 1902; Autonne 1904, 1905; Séguier 1902, 1907; Ranum 1907).

---

[35] Although Wedderburn was Scottish, he stayed at the University of Chicago at the beginning of the early century and worked in close collaboration with Dickson on the theory of algebras. It is in this context that Wedderburn proved that a finite field is necessary commutative. See [Parshall, 1985] and [Fenster 1998].

[36] The question of the articulation of this collective of texts with sociological or institutional identities is difficult and will therefore not be dealt with here. It would indeed require further investigations on algebra and number theory at the turn of the 20th century. Very little is known about the complex situation in France [Goldstein 1999]. Such issues are at the core of the collective ANR project CaaFÉ.

After 1910, the references to Jordan in papers classified in group theory fell down to 7%, i.e., to a weight similar to the one group theory already had in the Jordan corpus before 1895. Moreover, references to Jordan in English papers fell much below average between 1910 and 1929. Some authors, who had been initially interested in other types of references to Jordan than linear groups, went back to these references after having interacted with the Dickson network. Three main causes of this phenomenon may be listed.

First, most papers classified in group theory from 1910 to 1930 were of a sporadic form and appealed to earlier types of references to Jordan, i.e., to primitive groups (Manning 1908-1918), collineation groups (Blichfeldt 1905, 1907; Mitchell 1911, 1913; Burnside 1911), the finiteness theorem (Bieberbach 1911; Frobenius 1911), and to modular equations (Plemelj 1923). As for the Jordan-Hölder theorem, it was mostly referred to in treatises or synthetic presentations, thereby following the traditional *synthetic-Traité* type of reference (Frobenius 1916; Baudet 1918; Baumgartner 1921; Krull 1926; Onofri 1928).

Second, the attention of the main authors of the Dickson network shifted to different (even though very close) issues such as the invariants of quadratic forms and their geometrical interpretations or associative algebras and continuous groups.

Third, new monographs, such as [Burnside, 1897], [Dickson, 1901] or [Séguier, 1904] could now be substituted to direct references to Jordan's *Traité*.

The specific identity of the Dickson network did nevertheless not vanish all the sudden in 1910. The expression Galois field was actually more and more widely used in the U.S.A. and would continue to be sporadically associated to Jordan's linear groups by new authors (Study 1908; Scaris 1911; Ostrowski 1914). Moreover, some practices inherited from Jordan's works continued to be used even though they were more rarely identified explicitly as such. Actually, the extensive use some Americans (and, to a lesser extent, some French) authors had made of the generic form of reference to *Jordan's linear group* did not circulate at a broader scale than the Dickson network. This generic form of reference therefore faded away after 1910. As shall be seen in greater details later, the legacy of the Dickson network would nevertheless play an important role in the universalization of the generic reference to the Jordan canonical form theorem in the 1930s. Some echoes of this legacy would still resound in the second half of the 20th century, such as in the longstanding use of the expression "Champ de Galois" in France in parallel to the expression "Corps fini."

### 4. On the relation Jordan-Galois (1870-1914)

We shall now change perspective by questioning the legacy of Jordan's presentation of Galois theory in the time-period 1870-1914. Most references to Galois in the corpus we have considered until now were related to Galois fields but not to Galois theory. For the purpose of investigating this situation further, I am now comparing the Jordan corpus to the 315 reviews that referred to "Galois" during the time-period 1870-1914 (the Galois corpus for short).[37]

---

[37] Unlike searches on "Jordan," the name "Galois" does not raise issues relative to synonyms and therefore allows automatic searches in the long run. Moreover, no more than a small number of about ten reviews referred to Galois

The diagrams n°6 highlight that even though the relative weights of the Galois corpus and of the Jordan corpus both vary between 1 $^0/_{00}$ and 7 $^0/_{00}$ of the *Jahrbuch,* their evolutions are quite different. Actually, apart from the Dickson network, the intersections between the two corpora are limited to a very few reviews. It was thus obviously not at all common to refer to Jordan when mentioning Galois in the 1870s and 1880s. Moreover, most of the rare reviews that referred to Galois in connection to Jordan either pointed to the *synthetic-Traité* type of reference (e.g., Janni 1874; Netto 1882; Pellet 1887; Bolza 1890, 1893) or to the *special equations* type (Kronecker 1875; Gierster 1881). In contrast, these reviews did not highlight the issue of the association of groups to equations. Jordan was actually neither mentioned in connection to Klein's approach to Galois groups nor in the review on Bachmann's 1881 presentation of Dedekind's Galois theory. In sum, none of the texts reviewed in either the Jordan corpus or the Galois corpus actually followed the *Traité*'s specific approach to the "*Irrationnelles*".

The diagram n°6 points to three specific time-periods in the Galois corpus between 1870 and 1940: the beginning of the 1880s, the turn of the century, and the late 1920s.[38]

### 4.1. An overview of the Galois corpus from 1870 to 1890

Because Jordan's *Traité* has often been considered as having unfolded the true content of Galois's works, the references to the latter have usually been assumed to have been more and more concerned with groups after 1870 in contrast with more classical concerns for equations. This presentation raises several difficulties. First, it must be pointed out that the number of references to Galois remained very small after the publication of Jordan's *Traité* (15 entries during the ten-year period 1870-1880), and, while they increased at the beginning of the 1880s (25 entries during the five-year period 1881-1885), a significant quantitative evolution would not occur until the mid 1890s (67 entries during the five-year period 1896-1900). Second, only a single review mentioned Galois in the chapter equation (Frobenius 1872, on issues related to monodromy and differential equations) before the works of Klein and Brioschi on modular equations in the late 1870s. Third, the main sections of the mathematical classification associated to the Galois corpus remained significantly varied during the whole time-period. Most of the reviews were associated to five sections: "History and philosophy", "Algebra" (1. Equations, 2. Substitutions), "Function theory" (2. special functions), and "Differential and integral calculus" (5. General differential equations).

- **1870-1885: equations, substitutions, and special functions**

An overview on the great variety of topics of the papers reviewed in the chapter "equations" makes it clear that Jordan's approach to Galois theory was unlikely to cause any dramatic evolution. As a matter of fact, the main authors of this chapter in the 1870s were Clebsch and Brioschi (invariant theory), Laguerre (successive approximations), Cayley and Klein, while

---

because of the retrospective addition of the classification AMS 2000 (i.e. Galois Theory). Interestingly enough, most of these texts were authored by Dedekind who works otherwise appears only once in the corpus.
For an investigation of the Galois corpus in connection to the development of Galois theory, see [Ehrhardt 2007].
[38] One must be cautious with the *Jahrbuch* database in the 1920s-1930s because, on the hand, of the delayed reviews on many papers that were published during World War I, and, on the other hand, because of the political situation in Germany in the 1930s [Siegmund-Schultze 1993].

the reviews' main references pointed to the names of Charles Sturm (Sturm's sequences), Hermite and Kronecker (algebraic forms and special equations), Abel, Laguerre, Gordan, Cauchy, and Lagrange.

In his review of Jordan's 1869 "Commentaires sur Galois," Netto described the "well known results of Galois" as concerned with "general equations" as opposed to both "special equations" and "numerical equations." Only a very few reviews of the Galois corpus were nevertheless concerned with general equations. These usually did not point to group theory but to the traditional notion of "Galois's resolvent" in the legacy of Lagrange's approach to equations. [39] But most reviews concerned with equations in the Galois corpus did not refer to general equations but to special ones, i.e., the modular equations of degree 5, 7, and 11. In this context, Galois's name was mainly used in relation to "Galois-Betti-Hermite's results" on modular equations (Jordan 1868; Krause 1873; Briot et Bouquet 1875; Klein 1877-1884; Brioschi 1878; König, 1879; Krey 1880; Dyck 1881).

- **The 1880s: the Galois groups**

In the 1880s, the two types of references mentioned above were fading away. In parallel, more and more reviews were referring to *the* "Galois groups" to designate the three groups associated to the three modular equations of degree 5, 7, and 11. This type of reference to Galois increased especially from 1879 to 1883 in connection to the works of Klein and his followers, as is exemplified by [Gierster 1881].

We have already seen that the special equations of elliptic functions had underlying them a transversal framework at a large European-scale. As a result, most of the entries of the Galois corpus from 1870 to 1885 circulated between the chapters "equations," "substitutions," and "special functions." Let us consider for instance [Klein 1879a] and [Klein 1879b], two texts published by the same author on the same topic, the same year, and in the same journal. While the former was classified by Netto in the chapter "equation," the latter was reviewed by Müller in the chapter "special functions." The 1889 discussion between Brioschi and Halphen on the division by 7 of the periods of elliptic functions provides another example of the transversal dimension of the special equations of elliptic/hyperelliptic functions in regard with the *Jahrbuch* classification. While [Halphen, 1889] was classified in the chapter "substitutions" by Netto, Broschi's response was classified in the chapter "special functions" by Müller.

The increasing references to Galois in connection to Klein may have caused echoes and reactions in other contexts. In 1881, Bachmann published an approach to Galois theory based on the Dedekind's lectures. In 1883, Picard published a note in which he proposed to develop for differential equations a "method" analogous to the "method of Galois" for algebraic equations, and therefore to state the analogous of the "fundamental theorem of Galois" for Lie's continuous groups of transformations [Archibald 2011].

---

[39] A Galois resolvent is a function of $n$ variables that takes $n\,!$ values under the action of $Sym(n)$. In the sense of Lagrange's 1770 notion of similar function, a Galois resolvent is similar to any rational function of the roots of an equation of degree $n$ and of its coefficients (recall that the roots are supposed to be distinct). For more details on resolvents in Lagrange's tradition, see [Tignol 2001]

One year earlier, Kronecker had discussed the relevance of the notion of "Galois groups" in regard with his own notion of "equation with affect." From this point on, Galois's name circulated within a network of texts almost systematically classified by Netto in the chapter "equations." Combined with the decline of the chapter "special functions," this situation was the main cause for the temporary importance taken on by the chapter "equations" within the Galois corpus.

Amongst the specificities of the latter chapter were, on the one hand, its presence throughout the whole period and, on the other hand, the major role it could potentially play by gathering more than 60% of the entries at certain periods of time such as from 1885 to 1890. After 1890, the relative importance of this chapter nevertheless fell back to a level of about 20% of the corpus.

### 4.2. Intertextual networks

As has been seen above, the *Jahrbuch* classifications are not providing a partition subtle enough for characterizing the evolutions of the Galois corpus from 1870 to 1890. A closer look at the intertextual references during this time-period highlights the presence of two main networks of texts centered on the works of two individuals: Klein and Kronecker. Each of these two groups was mostly active in a single ten-year period even though some papers would continue to be published sporadically later on.

On the one hand, Klein's network was mostly active from 1877 to 1885 with a peak of activity between 1878 and 1882. It revolved around Klein's 1884 *Vorlesungen über das Ikosaeder und die Auflösung der Gleichungen vom fünften Grade.* On the other hand, the Kronecker network was especially active from 1880 to 1890 with a peak between 1882 and 1885. It revolved around Kronecker's 1882 *Grundzüge einer arithmetischen Theorie der algebraischen Grössen,* an alternative theory to Dedekind's approach to Körper and Ideals.

The few references to Galois that were independent from both the Kronecker and the Klein networks were also independent from Jordan's approach. All of them made long-term references to Galois's original works. Hermite's 1885 paper on Galois's continuous fractions was explicitly related to the curriculum of the teaching of mathematics [Goldstein, 2011]. The rare papers that still commented on Galois's criterion of solvability of irreducible equations of prime degree (Paxton Young 1885; Dolbnia 1887; Tognoli 1890) had underlying them the broad public dimension of the Galois criterion since Liouville's 1846 edition of Galois's works.[40] Until the mid 1890s, it was indeed mainly in connection to this criterion that Galois was publicly referred to in publications targeting larger audiences than specialized mathematical journals [Brechenmacher 2012].

- **The Klein network**

Most of the papers of the Klein network were published in German in *Mathematische Annalen,* i.e., Klein's journal. The main authors were Klein, Brioschi, Bachmann, Weber, and

---

[40] The concluding theorem of Galois's *Mémoire* is the criterion that: "in order that an equation of prime degree be solvable by radicals, it is necessary and sufficient that, if two of its roots are known, the others can be expressed rationally" [Galois 1846, p. 432].

Heinrich Maschke, as well as most of Klein's doctoral students in the early 1880s, e.g. Dyck, Gierster, Ernst Fiedler, Georg Friedrich, Robert Fricke, Adolf Hurwitz, and Frank N. Cole. But the Klein network also included more unexpected authors such as Georges Halphen and Poincaré.

Here, the references to Galois were usually not directly connected to Jordan's presentation of Galois theory. We have seen that Jordan had inscribed his *Traité* in a broad collective dimension by appealing to both the legacies of Hermite's and Clebsch's approaches to the special equations of elliptic/hyperelliptic functions. The references to Galois in Klein's network were of a similar type: they formulated Hermite's approach of the three modular equations of degrees 5, 7, 11 in the framework of some Clebsch-like geometric interpretation of the Galois resolvent, which Klein had actually opposed to Jordan's "abstract" approach to substitution groups in 1871 [Brechenmacher 2011a]. In this context, *the* Galois groups designated the groups of binary linear fractional substitutions associated to the three modular equations, i.e., $PSl_2(p)$ with $p=5,7,11$. Later echoes of this type of reference to Galois would resound at the turn of the century in some works on collineation groups in projective geometry.

As was said before, Klein's approach to binary or ternary linear groups was very different from Jordan's considerations on general linear groups. Moreover, in regard with modular equations, the *Traité* was only one reference among many others and was therefore usually not mentioned at all in the *Jahrbuch reviews.* On the other hand, Jordan's approach to substitutions groups was explicitly recognized as specific. But this specificity was not connected to Galois. Jordan's *Traité* was actually mostly referred to for some general notions, such as the ones of monodromy groups or isomorphisms. Amongst Klein's students, Dyck and Gierster studied Jordan's *Traité* in details along with the works of other authors such as Mathieu or Serret. As was said before, it was through Gierster that Jordan's linear groups would eventually circulate to the Dickson network at the turn of the century.

Even though they did not appeal to Jordan's presentation of Galois's theory, Klein and his followers nevertheless pursued a general aim similar to the one of Jordan's "*irrationnelles.*" Indeed, the Klein network revolved around the "fundamental problem" of the essence of the "irrational quantity" of the general quintic. The irrationality was represented by the "icosahedron," which could be considered as a polyhedron, a Riemann surface, an algebraic form, an equation, a transformation group, or a substitution group. One of the main motto of the Klein network was indeed to connect various parts of mathematics.

The Klein network was also interacting with other collectives of texts such as French works on linear differential equations (Poincaré, Picard, Halphen, Appel, Goursat, etc.). Here it is interesting to note that through these interactions, Klein's references to "Galois resolvents" or to the three "Galois groups" of modular equations circulated to [Poincaré 1883] and [Halphen 1889].

In his review of the first edition of Klein's 1884 *Vorlesungen*, Lampe especially highlighted the role the book devoted to "the theory of Galois groups." The designation still implicitly pointed to the three modular groups. But Klein's monograph also presented a general exposition of "Galois groups of equations" of degree $n$. However, this presentation was not based on Jordan's *Livre III* but on Kronecker's *Grundzüge*.

- **The Kronecker network**

The identity of the group of texts under consideration here is close to what one could designate as Kronecker's school in Berlin in the 1880s. Here, no reference to Galois was connected to Jordan's *Traité* even though Kronecker and his followers were mainly referring to Galois's theory of general equations.

Not only were the main authors of this network all former students of Kronecker, but most of the texts were published either in *Crelle*'s (i.e. Kronecker's) Journal or at the Academy of Berlin. The specific relation to Galois that circulated at the local level of the Kronecker network was indeed not shared at the larger scale of the *Grundzüge*'s influence. Making little reference to Galois and only to specific aspects of the latter's work was one of the characteristics of the group. Frobenius, for instance, would make almost no reference to Galois until the beginning of the 1890s. As for Netto, Adolf Kneser and Kurt Hensel, they would adopt their master's notions of Galois genus, Galois equation, and "equations with affects."

In the early 1850s, Kronecker had already claimed that it was impossible to fathom the "true nature" of solvable equations "from Galois's investigations. For Galois only addresses the first task to find the 'conditions of solvability'" [Petri and Schappacher 2004, p. 233]. Kronecker considered the theory of equations as resorting to the investigation of the "true nature" of irrational quantities. One will thus investigate equations of a given degree on the model of the explicit expressions Abel had given to the roots of the quintic, i.e., in the aim of finding the "most general" function by which the roots of any equations of a given degree could be expressed.

Kronecker's approach was therefore not very compatible with Jordan's focus on non-effective procedures on classes of groups and on equations of arbitrary degree [Brechenmacher 2007]. Kronecker and Netto were nevertheless close readers of Jordan's *Traité.* But unlike Jordan's *Livre III,* Kronecker clearly separated the effective arithmetic foundational issues on algebraic quantities he dealt with in his 1882 *Grundzüge* from the substitutions to which Netto devoted a monograph in 1882. When they were dealing with substitutions, both Kronecker and Netto appealed to the traditional point of view of the problem of the number of values of functions, i.e, to substitutions acting on functions of $n$ variables.

The traditional dimension of Netto's perspective on Galois theory is well illustrated by the important role the latter attributed to the "metacyclic equations" Galois had considered when he had stated his criterion of solvability of irreducible equations of prime degree. In presenting such equations as generalisations of cyclic and abelian equations, Netto was following Serret's *Cours d'algèbre supérieure*. Actually, apart from Jordan's *Traité* and Klein's *Icosahedron,* the solvable prime degree "metacyclic equations" would conclude most presentations of Galois theory until the turn of the century. Moreover, Netto presented the notion of Galois group as a secondary notion as compared to Kronecker's "concrete" notions of "Galois equation" (i.e., an effective reformulation of the notion of Galois resolvent) and of "equations with affects" (i.e., equations whose resolvents split up arithmetically in irreductible factors in a given "rational domain").[41]

We have seen that, unlike Kronecker, Klein had attributed a fundamental status to the Galois groups. The extensive and various uses of Galois's name that developed in connection to Klein's work circulated on a larger-scale than the Klein network (Galois ideas, Galois theorem, Galois groups etc.). Until the 1890s, most presentations of Galois's theory of equations would be based on Kronecker's approach. They would nevertheless follow Klein in attributing to Galois groups a status Kronecker had denied them.

- **Beyond Klein and Kronecker**

In sum, the circulation of Galois's theory of general equations was neither related in an obvious way to the reception of the *Traité* nor to group theory or even to algebra. Actually, auhors who, like Netto, Kronecker, and Pellet (1879, 1889, 1891),[42] were interested in Galois theory rejected Jordan's approach and carefully separated groups from the arithmetic dimension underlying "Galois's algebraic principles." But bits and pieces of Jordan's approach to the "*irrationnelles*" given by special equations nevertheless played a mediating role between Klein's and Kronecker's approaches. For instance, when Weber investigated some issues related to the double points of an algebraic curve (a quartic) in 1883, he presented his paper as concerning the "Galois group of an equation of the 28th degree."

Later on in 1889, Hölder appealed to Jordan's *Traité* to propose a synthesis of the works of Klein and Kronecker through what is now designated as the Jordan-Hölder theorem.[43] Hölder initially aimed at providing a new proof of "Galois's fundamental theorem" on the association of groups to equations, in response to a recent incorrect proof published by Söderberg (1887). His approach has often been celebrated for the new light it has shed on the chain reduction of groups by highlighting the fundamental role played by the concept of quotient group. Hölder's memoir was indeed announcing a larger-scale process of reorganisations and of institutionalization of group theory in the 1890s. In this context, new readings of Jordan's *Traité* developed in the lights of more recent achievements as has been discussed above in connection to the Dickson network. This process of institutionalisation especially gave the actors a free hand in dealing with the legacies of prominent authors of the 1880s such as Kronecker, Klein, and Lie.

In 1889, Hölder had divided his paper in two parts. The first was devoted to "pure group theory" in the legacies of Klein, Dyck and Gierster. The second was of an "algebraic" nature and presented Kronecker's approach to Galois theory. But what were exactly the interrelations between "pure group theory" and "algebra" ? As shall be seen in the next paragraph, at the turn of the century, the disciplinary issues related to the nature of Galois groups would often have national, and in fact nationalistic, overtones.

---

[41] In short, Kronecker had developed a constructive presentation of finite field extensions of certain ground fields. See [Petri and Schappacher 2007], [Goldstein and Schappacher 2007 p. 81-88]

[42] On Pellet's approach to Galois's work, see [Ehrhardt 2007].

[43] For a comparison between Jordan and Hölder's approaches on Jordan-Hölder theorem, see [Nicholson 1993] and [Corry, 1996, p. 24-34].

### 4.3. Disciplines (1890-1914)

The main cause for the relative decline of both the chapters "special functions" and "equations" at the turn of the century was the parallel growth of four other chapters: group theory, number theory, differential equations, and history (diagrams n°7). We have seen in the previous section that most types of references to Galois were mainly transversal to the *Jahrbuch* classification in the 1870s and 1880s. At the turn of the century however, each of the three disciplines of algebra, arithmetic, and analysis, pointed to specific types of references to the relation Jordan-Galois.

- **Algebra**

Although the algebraic "group theory" dominated the whole corpus, this domination was mainly due to the Dickson network (85% of the papers in "group theory"). As seen above, linear groups in Galois fields pointed directly to Jordan's *Traité,* i.e., neither to Kronecker's general theory of equations, nor to the three Galois groups of Klein's network.

- **Number theory**

References to Galois in "number theory" increased after the publication of Weber and Hilbert's 1893-1894 "new grounds for Galois theory." It is not the place here to stress an overview of the development of algebraic number theory in Germany.[44] We shall nevertheless illustrate the increasing algebraic-arithmetical disciplinary framework developped in the texts in German language of the Galois corpus by the evolution of three successive synthetic works published by Weber.

The 1885-1887 "Zur Theorie der elliptischen Functionen" gave of presention of both the Galois groups of special equations and the Galois principles for general equations. But these two lengthy papers published in *Acta mathematica* still appealed to a transversal approach in the legacy of Kronecker's complex multiplication of elliptic functions and its relations to quadratic forms and class-numbers. In contrast, Weber's 1891 *Elliptische Functionen und algebraische Zahlen* aimed at "giving for the first time a comprehensive development of the relations between elliptic functions, algebra, and number theory." The book was therefore partitioned into three sections, each devoted to a discipline. References to Galois nevertheless still played a transversal role in each section.

The first volume of the well known *Lehrbuch der Algebra* appeared in 1895. Even though it was implicitly still resorting to the legacies of Klein and Kronecker, this book emphasized Dedekind's approach to Galois theory, i.e., the interplay of extensions of "Körpers" and of the resulting "division" of a groups into normal subgroups. In 1897, Hilbert's *Zahlbericht* presented for the recently created *Deutsche Mathematiker-Verinigung* the arithmetic theory of a general Zahlkörper *K* with a detailed investigation on the decomposition of the prime ideals of *K* in a Galois extension of *K*. At the turn of the century, Galois theory was thus eventually considered as an elementary part of "algebraic number theory" which identitified itself to both a discipline – Algebra, with a clear arithmetical trend -, and a national frame.

---

[44] Cf. [Corry, 1996], [Goldstein and Schappacher 2007b], [Petri and Schappacher 2007].

Almost all the texts of the Galois corpus that appeal to the notion of "Körper" were published in Germany. Most of them were classified in the section "number theory" of the *Jahrbuch.* They actually amounted to more than 75% of the texts of the Galois corpus classified in the chapter "number theory. 1. Generalities." The authors appearing in this context are Weber Hilbert, Minkowski, Bachmann, Dedekind, Hensel, Sapolski, Bauer, Fueter, Frobenius and Furtwängler. At the dawn of World War I, the traditional sections "Algebra" and "Niedere und höhere Arithmetik" of the classification of the *Jahrbuch* were merged into a section "Arithmetik und Algebra" with a subsection devoted to Galois's theory.

- **Analysis**

We have already seen that, except in the case of Galois fields, the circulation of references to Galois was neither related in an obvious way to the reception of Jordan's *Traité*, nor to the development of group theory. But we have seen also that the reception of Jordan's *Traité* was not limited to algebra or group theory. In France, this reception had indeed especially involved linear differential equations and the theory of algebraic forms. But in neither of these two contexts did the authors refer to Galois when working with Jordan's method of substitutions. Let us consider the example of Poincaré's works in the 1880s. Although the appropriation of Jordan's method had played a key role in the development of the theory of Fuchsian functions [Brechenmacher 2011b], the rare occasions when Poincaré referred to "Gallois" (sic.) were actually related to the works of Klein.

Moreover, apart from the context of the teaching of the *Algèbre supérieure*, there was at the time no shared reference to Galois's works in France. This situation is illustrated by Léon Autonne's unique reference to Galois in 1885. Even though Autonne had been initiated to Jordan's substitutions by Poincaré, Autonne never used any of the types of references to Galois that were circulating in Klein's network (to which he was not connected). The expression "Galois equations" he appealed to in 1885 was taken directly from Jordan's *Traité*. It shows that Autonne had not read Netto's 1882 treatise on substitutions. The latter had indeed replaced Jordan's "Galois equation" by the expression "metacyclic equation" because, in Kronecker's terminology, a Galois equation actually designated the resolvent of an equation.

Amongst the French mathematicians who were working on linear differential equations, Picard was the only one to develop a specific way to refer to Galois. Just after he had met Lie in Paris, he claimed in 1883 to develop for differential equations an approach "analogous" to Galois theory of general equations. References to a method "analogous" to "Galois's method" in the case of differential equations would developped in the mid 1890s in a group of texts involving authors such as Picard, Vessiot, Drach, Guldberg, Beke, Marotte, Cartan,

Schlesinger, and Loewy.[45] In this context, Picard's recurrent references to the "Galois's ideas" in relations to the notions of "groups" and of "general irrationalities" indicate the unifying role the latter attributed to Analysis.

At the turn of the century, several authorities such as Picard, Jules Tannery, Poincaré, and Jacques Hadamard contrasted the "richness" of the power of unification of Analysis with the "poverty" of considering algebra and/or arithmetic as autonomous disciplines. These official lines of discourse usually pointed to recent developments in Germany in the legacies of Kronecker or Richard Dedekind. In France, Galois was then increasingly celebrated as the founder of the analytic notion of continuous group.

Recall that for most of the 19[th] century in France, the theory of equations was no more an autonomous domain of research than algebra itself was an object-oriented discipline shared by a community of specialists [Brechenmacher and Ehrhardt 2010]. The *Algèbre supérieure* had developped in the 1840s as an intermediate discipline between elementary arithmetic/algebra and the "higher" domain of analysis as it was taught at the *École polytechnique* [Ehrhardt 2007, p. 211-236]. At the turn of the century, while Drach's *Algèbre supérieure* had appealed to Kronecker for presenting the "famous theory created by Galois" as an extension of arithmetic, Tannery's *Préface* explicitly recalled that it was Analysis that provided a higher point of view on "the general irrationality" of which "the algebraic number is nothing more than a particular case" [Borel et Drach 1895 p. iv]. The example of Picard, who included a presentation of Galois Theory in the third volume of his *Traité d'analyse*, makes it clear that the algebraic presentation of Galois theory was considered as a first step toward the higher point of view of analysis. Even later, in 1913, George Humbert's lectures on the theory of substitutions at the *Collège de France* focused on the notion of group of monodromy and insisted along the lines of Jordan's *Traité* on applications in analysis and geometry.

Recall that Jordan's *Traité* had already presented Galois theory as a general theory aiming at providing a higher point of view on the "irrationals." The general part of this theory was explicitly considered as belonging to the Analysis, it was then to be applied to algebra, geometry and transcendental functions. It now appears that this presentation was very coherent with dominant views in France on the essential and transversal role of analysis on the disciplinary organization of the mathematical sciences.

- **History**

The boom of the section "History" was another important evolution of the Galois corpus at the turn of the century. This evolution points to the larger-scale phenomenon of the increasing relative weight of the section "History" in the *Jahrbuch*, i.e., from about 7% in 1871-1875 to a about 13% in 1911-1914. Complementary to the increasing number of academic publications (obituaries, collected works etc.), the faster pace of growth of the section "history" - in regard with the one of the whole *Jahrbuch*  - was also due to the development of histories of "disciplines." Moreover, the turn of the century was also a period of growth of the number of students of mathematics; a number of publications on the history of mathematics were participating to the construction of the curriculums of academic disciplines.

In this context, Galois was increasingly considered as one of the main founders of group theory [Ehrhardt 2007, p.628-649]. But the Galois corpus also highlights the various roles taken on by the history of mathematics at the turn of the century. In 1895, Sophus Lie had been invited to celebrate "Galois's influence on mathematics" at the occasion of the centenary of the *École normale supérieure*. Two years later, the *Societé mathématique de France* had Galois's works reprinted. Picard's introduction followed the role Lie had attributed to Jordan, namely as the one who had clarified and generalized Galois's distinction between simple and compound groups to the notion of composition series. Picard's claims were to circulate at an international level. Altogether with Klein, Picard and Lie played a key role in the consideration of Galois as the main founder of group theory in regard with other mathematicians such as Lagrange or Cauchy. After 1897, most histories of the theory of equations would usually adopt a three-act structure: before Galois, Galois, and how Jordan had made Galois theory "become public" [Pierpont 1897, p.340].

But the role of the researcher who closed the algebraic issue of the solvability of equations was also, somewhat incidentally, assigned to Jordan. Lie and Picard indeed both claimed that Galois groups had exceeded the boundaries of algebra. This role that was attributed to Jordan at the turn of the century can thus hardly be disconnected from the role of public authority on mathematic Picard had taken on at the turn of the century. Both the role assigned to Jordan and the one taken on by Picard were indeed fitting into a type of public claims on the unifying power of analysis which most authorities of French mathematics emphasized from 1880 to 1930. Even though such claims would involve some nationalistic discourse, in which Galois, like other *grands savant,* was much involved, they were not limited to anti-German discourses but pointed also to the constitution of a public expression of mathematics in France [Brechenmacher 2012].

Despite these public claims on the relation Jordan-Galois, we have seen that the increasing importance attributed to Galois's works in the late 19[th] century was not directly related to the reception of Jordan's *Traité*.

### 5. Quantitative perspectives on the time-period 1910-1939

Let us now return to the Jordan corpus for investigating further the time-period 1910-1939. The crystallization of some generic references to Jordan at the turn of the century would have a major impact on the Jordan corpus after 1910. The distribution of the languages in the Jordan corpus from 1910 to 1930 reflects the one of the references to "Jordan's curve" in the chapter "Topology": French and German languages are overrepresented. Papers in English language are catching up in the mid 1920s and the distribution of languages in the 1930s is close to the one of the whole *Jahrbuch* (diagrams n° 4).

The time-period 1910-1939 witnesses a reduction of most references to Jordan to a few generic designations such as Jordan's curves, Jordan's measure, Jordan-Hölder, and Jordan's canonical form. In contrast with this reduction, such expressions were used in an increasing diversity of chapters and languages. Indeed, not only were the five main sections of pure mathematics represented in the Jordan corpus but also more than 20 subchapters from 1910 to 1929 and more than 30 subchapters from 1930 to 1939. This situation is not surprising. Recall that the complex of notions at the roots of the Jordan curve theorem (such as laces of integration in complex analysis) were already thought to be transversal to several disciplines in the 1860s, especially in connection to the *regular solids* research field.

### 5.1. Jordan's

Let us take a closer look at the generic references to the Jordan curve theorem. In a series of papers he published between 1896 and 1904 in the *Göttingen Nachrichten,* Arthur Schönflies repeatedly referred to the "Jordan theorem" on the conditions for a curve to divide the points in space into two categories depending it they are inside or outside the curve (Schönflies, 1896, 1899, 1902, 1904).[46] Schönflies' works were stemming from the specific mix of concerns that we have seen to be characteristic of the *regular solids* research field, i.e,. motions, polyhedrons, symmetries, crystallography, and surfaces [Schönflies 1891a,b, 1892, 1893, 1896]. They triggered further discussions on "Jordan's curve" in the *Göttingen Nachrichten* (Bernstein 1900, Osgood 1900, Hilbert 1901).

The concept was shortly afterward used by Osgood (1902) in the *Transactions of the American Mathematical Society,* and it would from then on be discussed quite often in the *A.M.S.* (Bliss 1903, 1905; Ames 1903, 1905; Veblen 1905). After Hilbert's "Ueber die grundlagen der Geometrie" had originally been published in the *Göttingen Nachrichten,* it was published again in *Mathematische Annalen* in 1902. The generic reference to the Jordan curve theorem then circulated to the latter journal (Riesz 1904). This notion was then involved in the definition of a line and a surface in the German *Enzyklopädie* der *Mathematischen* Wissenschaften (Mangoldt 1905) and was quickly commented on by papers published in Italian and in French (Sbirani 1905; Zoretti 1907). The uses of this generic reference would quickly increase after the new developments on "Jordan's theorem" by actors such as Brouwer (1910, 1912), Lebesgues (1911), or Denjoy (1911) in the context of the development of topology as a mathematical discipline.

Let us not come to some more general conclusions about the introduction of generic references. First, it must be pointed out that the designation of an object by Jordan's name always resulted from a major evolution in the collective organization(s) of both the texts and the themes to which this object had used to be related. In a word, generic references attributed a collective meaning to something that did not make sense collectively any more.

For instance, the fast circulation of Schönflies's reference to the Jordan curve theorem in the context of discussions on the foundations of geometry was contemporary to the decline, if not the disappearance, of both the *regular solids* research field and of Klein's approach to *special equations*. While in these earlier collective organizations of knowledge, precise aspects of Jordan's works used to be discussed in connection to some specific texts by other authors, the objectified Jordan's curve, however, would be discussed in some new contexts, such as the definition of the notion of a set of points in connection to the works of Hankel, Cantor, Peano, and Borel (Schönflies 1900, Moore 1901, Gundersen 1901, Vitali 1904, Sbirani 1905).

---

[46] See [Scholz 1989, p. 120-124 & 137-148].

The Dickson network provides another example of this situation. Here, a collective of texts actually crystallized itself in attributing to Jordan's linear groups a generic nature after having used to refer individually to some specific parts of Jordan's *Traité*. At the turn of the century, several references to Jordan were indeed still presenting a mixed nature. In group theory, while most authors who were referring to the Jordan-Hölder theorem were not appealing to Jordan's works any longer, other references were still pointing directly to specific parts of Jordan's works, such as to the finiteness theorems, periodic substitutions, etc. Similarly, while authors such as Baire (1899) and Lebesgues (1904) were discussing the "Jordan measure" in connection to Borel's works, other authors were still referring directly to the *Cours d'analyse* (Böcher 1895; Pringsheim 1898; Moore 1901; Stolz 1902; Hardy 1903; Richardson 1906; Hobson 1906; Porter 1907).

The objects identified by generic designations were much more malleable than specific references to Jordan's original writings. The Jordan-Hölder theorem exemplifies how such objects could evolve and even take on various meanings simultaneously in various contexts. Most references to this theorem from 1910 to 1930 indeed occurred outside of group theory. The generic reference to Jordan-Hölder indeed carried on analogies that extended the process of chain reduction of a group to other objects, such as associative algebras (Epsteen and Wedderburn 1905), linear differential equations (Loewy 1912), or ideals (Sono 1924). Later on, the Jordan-Hölder theorem would play a model role in *Moderne Algebra* (Noether 1926; Krull 1928; Schmidt 1928; Schreier 1928). As was already the case with Jordan's curve theorem, echoes of the transversal nature Jordan had originally attributed to his method of reduction were thus resounding in the 1920s. After 1930, the Jordan-Hölder theorem represented more than 85% of the references to Jordan in the chapter "Group theory / abstract algebra" while only five papers kept making sporadic references to Jordan's works, i.e., to the limit of transitivity (Weiss 1930), the origin of the linear group (Bottema 1930), the synthetic nature of the *Traité* (Got 1933), and the finiteness theorems (Zassenhaus 1938; Turing 1938).

### 5.2. The universalization of the Jordan canonical form theorem

Both the Jordan curve theorem and the Jordan-Hölder theorem point to a phenomenon of reduction to a generic designation of a variety of references to some specific parts of Jordan's works, combined with an increasing diversity in topics and languages. I shall designate here such a phenomenon as a process of universalization. As is illustrated with the examples of the generic references to "Jordan's method" in the 1880s or to "Jordan's linear groups" at the turn of the century, the development of a generic reference does not imply directly any process of universalization.

This section aims at investigating further the process of universalization of the generic designation to "Jordan's canonical (or normal) form" theorem. Unlike the Jordan curve theorem, the references to the latter theorem were mostly sporadic until the 1930s even though they had turned into some continuous, and even generic, forms of references in some local contexts. Recall that designations such as Jordan's curves and Jordan's measures had been coined in Göttingen, and in Paris, two of the main centers of the production of mathematics at the time. In contrast, we have seen that the circulation of the expressions Galois fields or Jordan's canonical forms was at first limited to the Dickson network.

In the 1930s, however, the expression Jordan's canonical form theorem was not only circulating at an international level, it was also invading a great number of sections of the *Jahrbuch* classification. This phenomenon must be considered in connection to the larger-scale process of the universalization of matrices and of the emergence of linear algebra as a discipline (a section "linear algebra" was created in the *Jahrbuch* in 1939 with a subsection devoted to "matrices"). This evolution resulted from a complex of factors. One may especially point to the combination of two phenomena. The first is the decline of the theory of bilinear forms after World War I, while since Frobenius's synthesis in 1877-1879 this theory had played a major role in algebra and had been applied to various topics (e.g. mathematical methods in physics, differential equations, integral equations, relativity theory, etc.). In the early 1920s, it was still one of the main elementary topics of treatises of algebra and was also underlying various researches. The second factor is the key role that was played by matrix mechanics in the quantum theories of the mid-1920s. Quantum theories indeed supported both the internationalization and the diversification of the uses of the terminology matrix in regard with other terminologies, such as the one of "*Tableau*" in France [Brechenmacher 2010].

But the universalization of Jordan's canonical form theorem points to a third factor, which is the globalization of some local specific uses of matrices or *Tableaux* through the internationalization of the theory of canonical matrices in the 1930s.

Let us first get an overview on the circulation of Jordan's canonical form from the 1870s to the 1930s. Jordan introduced his canonical form for linear substitutions on integers mod. $p$ (or mod. $p^n$) between 1868 and 1870 [Brechenmacher 2006]. The theorem embodied the method of reduction we have seen to be specific to Jordan. It especially resorted to the unscrewing into the two forms of actions of cycles (*i gi*) and *(i i+a)* which it assimilated to issues involving $n$ variables. Later on, Jordan would appeal frequently to the canonical reduction of substitutions (in *GF(p^n)* or $\mathbb{C}$) in his works on groups, differential equations, algebraic forms, etc.

But the theorem nevertheless almost disappeared from the public scene after it had been strongly criticized by Kronecker in 1874 [Brechenmacher 2007]. Kronecker not only rejected the formal generality of Jordan's linear groups, but he also criticized the non-effectiveness of the canonical reduction, which required the determination of the roots of arbitrary algebraic equations. Unlike Jordan's theorem, Kronecker promoted the computations of invariants by determinants such as Weierstrass' 1868 elementary divisors [Hawkins 1977], which he reformulated as a rational method of computations of invariants in any "domain of rationality" (i.e. the invariant factors of matrices in a principal ideal domain).

Recall that various ways of dealing with linear substitutions had parallel circulations until the constitution of linear algebra as a discipline in the 1930s [Brechenmacher 2010]. Amongst these, the most influential approach was based on Frobenius's 1877-1879 presentation of bilinear and quadratic forms. This approach appealed to symbolic methods and to computations of invariants by determinants [Hawkins 2008]. But Frobenius had presented Jordan's canonical form as a corollary of Weierstrass' elementary divisor theorem. Following Kronecker, he moreover insisted that the validity of Jordan's form was limited to the case when one would allow the use of "irrationals" such as "Galois's imaginary numbers" [Frobenius 1879, p. 544].

As a result, the circulation of Jordan's canonical form went underground in the works of authors such as Poincaré or Élie Cartan, where it was neither considered as a theorem nor attributed to Jordan [Brechenmacher 2012]. From 1875 to 1895, the canonical form theorem can thus hardly be detected in the Jordan corpus. Gaston Darboux referred it to explicitly in 1874 in the context of the Jordan-Kronecker controversy. An implicit reference was also made by Hamburger in 1881 in a review of a paper of Casorati on linear differential equations. Hamburger indeed noted that the author had developed a different approach than Jordan's by appealing to Weierstrass's elementary divisors. [47] In 1889, Predella mentioned that even though he had appealed to Jordan's $n$-dimensional geometric space, he had nevertheless developed an alternative method by appealing to Weierstrass's theorem.

At the turn of the century however, Jordan's canonical form was circulating in plain sight in the Dickson network. [48] Much work was devoted to making some procedures of matrix decomposition explicit that had never been considered as mathematical methods *per se* until then (Burnside 1898, 1899; Dickson 1900-1902: Séguier 1902, 1907; Autonne 1905-1910; Châtelet 1911). Unlike the static nature of the invariants of the Frobenius theory, Jordan's canonical form was indeed based on some dynamic decomposition of the analytic representations of matrices. As was claimed by Séguier at the Academy of Paris in 1907, and later on by Dickson in 1924 at the congress of mathematicians in Toronto, attributing a key role to Jordan's canonical form implied reorganizing the structure of Frobenius theory of bilinear forms in laying the emphasis on transformations of matrices rather than on classes of equivalences of forms [Brechenmacher 2010].

After 1910, Jordan's canonical form circulated in an international network of texts in which the "theory of canonical matrices" was arising as an autonomous topic of research. Concerns for the teaching of mathematics were at the origin of the development of this network. Indeed, most of its texts referred to a paper published by Lattès in 1914. The latter nevertheless did not refer to the Dickson network but to the use Jordan had made of his canonical form in his *Cours d'analyse*. This situation illustrates the mixed status of Jordan's canonical form at the beginning of the 20[th] century. The theorem was indeed at the same time a generic reference in the Dickson network, a continuous reference in the framework of the teaching of differential equations, and either a sporadic or non-existent reference in most other contexts.

---

[47] Hamburger had actually been the first to point out the relation between Jordan's canonical form and Weierstrass' elementary divisors in 1872 [Brechenmacher 2007].

[48] The somehow discontinuous circulation of Jordan's practice of reduction from the 1860s France to the 1890s Chicago is analyzed in [Brechenmacher 2011a].

Lattès's paper provided an iterative effective procedure for computing Jordan's canonical form. Its purpose was to set a method that may be used for the teaching of linear systems of differential equations with constant coefficients. But Lattès also incidentally tackled the main issue of the 1874 controversy between Jordan and Kronecker. He indeed showed that one could both adopt Kronecker's ideal of effectiveness and Jordan's practice of reduction of a substitution to its simplest form.[49]

That the theory of canonical matrices was both related with teaching issues and pointing to open research questions was a part of the process of universalization of matrices in the interwar-period. Indeed, complementary to the increasing uses of matrices in various topics and languages, matrices were also getting involved at the various levels of mathematical activities. Following Lattès, a series of texts investigated rational procedures of canonical decompositions (Kowalewski 1916; Loewy 1917; Krull 1921; Dickson 1926; Polya 1928; Bennett 1931; Rutherford 1932; Ingraham 1933; Mac Duffy 1933; Gantmacher 1935; Cavalluci 1937), while others focused on non-rational procedures of reductions to what they were designating as the "classical canonical form," or as "Jordan's canonical form" (Burgess 1916; Voghera 1928; Aitken 1928; Wellstein 1930; Bell 1930; Smale 1930; Turnbull 1931; Menge 1932; Amante 1933; Cherubino 1936; Cramlet 1938).

In the 1930s, these issues were interlaced to some other uses of matrices, for instance in connection to the arithmetic of associative algebras in principal ideal domains (Dickson 1926; Wedderburn 1931; Mac Duffee 1933; Ingraham & Wolf 1937), to the fundamental theorem for finitely generated abelian groups (Châtelet 1922, 1923; de Séguier 1925), or modules (Krull 1926; Van der Waerden 1931), or to some issues on group representation related to quantum mechanics (Weyl 1927; Wintner 1927; Fantappié 1928; Van der Waerden 1932; Bauer 1933; Albert 1934; Schwerdtfeger 1935).

These connections and reorganizations are exemplified by the number of distinct solutions that were given to the problem of the determination of the matrices that commute to a given matrix (Kravčuk 1924; Shoda 1929; Bell 1930; Rutherford 1932; Williamson 1934; Hopkins 1934). The latter issue had actually been tackled a great number of times since the 1850s in connection to various objects: substitutions, algebraic forms, matrices, associative and Lie algebras, groups, etc. Its presentation in the framework of the theory of matrices played an important role for the unification of the latter theory.

---

[49] In modern parlance, Lattes reduced a matrix in a principal ideal domain to a chain of companion matrix. The latter notion had been introduced by Frobenius in 1879. In Frobenius theory, companion matrices were nevertheless not deduced from some matrix transformations but from some computations of invariants (Kronecker's notion of invariant factors of a matrix in a principal ideal domain). Some rational methods of reduction to canonical forms similar to the iterative methods of Lattès had already been given prior to the latter by Landsberg and Burnside at the turn of the century. These methods were presented in a treatise published by Hilton at about the same time as Lattès paper in 1914. It was nevertheless eventually the latter that became the main reference in the development of the theory of canonical matrices.

In the 1930s, decompositions to canonical forms laid the ground for most presentations of the theory of matrices. In this context, the relations between Jordan's and Weierstrass' theorem were discussed again (Lusin 1931, Mac Duffee 1933), and even though the expression "Jordan-Weierstrass's theorem" was used from 1936 to 1938 (Julia 1936, Zwirner 1936, Rothe 1936; Volterra & Holstinsky 1938), the former eventually dethroned the latter.

This phenomenon sheds a new light on the modernity of algebra at the time. The historiography of algebra has indeed usually focused on the emergence and the diffusion of the *Moderne Algebra* through van der Waerden's presentation of Emmy Noether's and Emil Artin's lectures. But the section on "linear algebra" in van der Waerden's book was nevertheless showing a very traditional structuration. It indeed followed the structure of Frobenius's theory of bilinear forms in appealing to the symbolic representation of matrices. In this context, the tabular form of matrices was limited to some *a posteriori* static representation of sequences of polynomial invariants – such as Weierstrass's elementary divisors. According to some of the French Rockefeller fellows in Hambourg and Göttingen, E. Noether was actually explicitly against the procedures of matrix decomposition [Goldstein 2009, p.165]. But the abstract ideal Noether emphasized in this specific context was a traditional one in the legacy of Frobenius theory. The latter theory had indeed incorporated the notion of matrix in the late 1880s. But matrices were mostly assigned the secondary role of an abstract notion underlying the symbolic operations that could be made on the distinct central objects of Frobenius's theory, i.e., bilinear and quadratic forms, substitutions, and determinants. Later on, the *Moderne Algebra* would assign the abstract role to algebraic structures such as moduls and rings, thereby turning Frobenus's matrices into a superfluous notion.

It is interesting to note that Cyrus Colton Mac Duffee, a former student of Dickson, was following closely van der Waerden's book in his first treaty on *The theory of canonical matrices* (1933). This book was indeed mostly devoted to translate some traditional results on matrices by appealing to the abstract notions of "linear algebras" (i.e. associative algebras), moduls, and ideals. Mac Duffee thus followed the structure of Frobenius's traditional theory even though he made a generic reference to the Jordan canonical form theorem in Dickson's legacy. But in contrast with the dynamic procedures of matrix decompositions that had been circulating in the Dickson network, Mac Duffee stated the Jordan theorem by appealing to some static invariant computations based on Weierstrass's elementary divisors theorem. In his 1941 treaty *Vectors and matrices,* Mac Duffee eventually abandoned van der Waerden's abstract approach for focusing on the complex dynamic procedures of transformations of the tabular forms of matrices and on their geometric interpretations as the action of an operator on a finite dimensional a vector space.

The universalization of the generic reference to the Jordan canonical form theorem thus sheds new light on some evolutions of algebra in the 1930s which have often been underestimated in regard with Göttingen's *Moderne Algebra*. Moreover, the connections between these various evolutions raise some complex issues. It is indeed not possible to oppose the modern abstract approaches of the ones to the concrete applications of the others. For instance, the procedures of matrix decomposition provided a model for the decomposition of a linear space under the action of an operator. The practices of matrix decomposition were thus already sewing the seeds for the ulterior decline of the matrix representation with the growing concerns for infinite dimensional spaces.

## Conclusion

In this paper, we have investigated the evolutions of some collective forms of references to Jordan from 1860 to 1940. We have appealed, on the one hand, to a global corpus of *Jahrbuch* reviews and, on the one hand, to some finer investigations of a few more local networks of texts. The *Jahrbuch* reviews thus appear to be an efficient tool for tackling some, otherwise unsolvable, issues such as the inverse intertextual relationships problem. They indeed provide a simplified overview through a quite homogeneous corpus. Although simplified, this overview gives access to a much more complex landscape than one may have initially suspected, especially given the individual nature the historiography has usually assigned to Jordan's works. But we have seen that the use of such a corpus shall not be limited to some quantitative analysis. It requires on the contrary a careful attention to the forms of distributions of textual references in time.

Among the various forms of distributions of references we have discussed in this paper, the sporadic form is especially important. We have seen for instance that the sporadic nature of some specific references to Jordan in a few reviews on the works of Poincaré, Picard, and Autonne has underlying it the main framework in which Jordan's linear groups were received in France, i.e., the context of linear differential equations. In contrast, this framework has remained invisible to the retrospective approaches that have looked for the reception of Jordan's *Traité* in some corpora that are fitting into nowadays algebra and finite group theory. It is important to point out that this framework is not directly visible in the Jordan corpus either: it was by following the lead shown by a sporadic form of references that this framework was identified through some finer investigations on the intertextual relationships of the papers under review.

We have seen also that most references that are contemporary to the papers they are pointing to take either a punctual or continuous forms. Continuous references are nevertheless doomed to turn either generic or sporadic at some point. This phenomenon requires further investigations. It is especially striking that no sporadic reference ever disappeared in the time-period the present case study has considered. For instance, we have seen that most early references to Jordan's *Traité* pointed to some *special equations* such as the equation to the 27 lines or the three modular equations Galois, Betti, and Hermite had especially investigated. We have seen also that these references were of a continuous form. They were indeed related to the "*irrationnelles.*" These involved especially elliptic and abelian functions and pointed to a collective dimension transversal to algebra, arithmetic, analysis, and geometry. In the 1870s-1890s, this transversal framework was torn apart by the development of some object-oriented disciplines such as finite group theory. The reference to Jordan's *special equations* then turned sporadic. But it never disappeared. This situation calls for some further investigations about the modes of transmissions of the memory of some mathematical issues. Sporadicity actually seems to be not only a form of distribution of references in time but also a form of transmission. Indeed, it seems that some specific results may be chosen sporadically as topics of investigations (especially doctoral investigations) precisely because of their sporadic nature: new light may indeed be shed from time to time on the works that are not fitting anymore into contemporary theories or disciplines.

The issue of the choice of a topic of investigation in connection with its distribution in time also calls for some further investigations on the *Jahrbuch* reviewers' modes of selecting and classifying papers. Reviews, selections, classifications are indeed forms of discourses on mathematics. They are instrumental in constituting a memory of some mathematical works and in transmitting such a memory. The population of the reviewers should especially be investigated further. Indeed, even though some individuals wrote occasional reviews, most papers were reviewed by a small number of actors (usually professors) who were specialized in a type of publications, i.e. usually a section of the classification, such as algebra for the case of Netto or special functions for the case of Müller. But some mathematicians may also have selected the papers they reviewed for other issues than thematic ones. Arthur Cayley, for instance, wrote hundreds of reviews of papers published in Great Britain. The main *Jahrbuch* reviewers wrote a large number of reviews a year (for instance, Lampe wrote more than 400 reviews in 1890). This impressive work gave them both responsibilities and power. For instance, in a letter he sent to Borchardt in 1874, Jordan complained about Netto's review on the *Traité*.

We have seen that some of the dynamics of the forms of references to Jordan's works are quite complex. The collective dimensions of these references contrast with the usual presentation of the relation Jordan-Galois as an exclusive one. This presentation has indeed supported some claims on the individual specificity, and in fact isolation, of Jordan's works [Julia 1962; Klein 1921], such as when Lebesgue considered that apart from his works on algebraic forms and substitutions groups, Jordan had tackled "weird" issues, "apparently disconnected one another" [Lebesgue 1923, p. 92]. We have nevertheless seen that issues such as $n^{th}$ dimensional geometry, crystallography, cinematic, topology, symmetries, etc., were fitting into a research field. We have seen also that the continuous forms of references to *regular solids*, *special equations*, or the *calcul des Tableaux* were transverse to disciplinary or national frameworks.

The local algebraic practices that have circulated over the course of the $19^{th}$ century in the "*calcul des Tableaux*" especially shed light on the reception of Jordan's *Traité* [Brechenmacher 2011a & 2012]. The globalization of such practices during the interwar period calls for further investigations on the generalization of the Jordan canonical form theorem in both the 1930s algebraic-arithmetic approaches to associative algebras, [50] and the algebraic-geometric-analytic approach to operator theory in (infinite dimensional) Banach spaces. As was already the case with Jordan's curve theorem, echoes of the transversal dimension Jordan had originally attributed to his method of reduction were thus resounding in the 1920s.

The transversal dimension of some generic references, such as the ones to the Jordan curve theorem or to the Jordan canonical form theorem, questions the non disciplinary patterns that mathematical disciplines have underlying them.

**Bibliography**

ADHEMARD (Robert d')
[1922], Nécrologie. Camille Jordan, *Revue générale des sciences pures et appliquées*, t. 3 (1922), p. 65-66.
ARCHIBALD (Thomas)
[2011] Differential equations and algebraic transcendents: French efforts at the creation of a Galois theory of differential equations (1880-1910), *Revue d'histoire des mathématiques,* à paraître.
AUTONNE (Léon)
[ ]
BACHMANN (Paul)
[1881] Über Galois' Theorie der algebraischen Gleichungen, *Mathematische Annalen*, 18 (1881), p. 469-468.
BERTRAND (Joseph)
[1867], *Rapport sur les progrès les plus récents de l'analyse mathématique*, Paris : Ministère de l'instruction publique, Imprimerie impériale, 1867.
BKOUCHE (Rudolf)
[1991] Variations autour de la réforme de 19021905, *in*  [Gispert 1991, p. 181-213].
BOREL (Émile), DRACH (Jules)
[1895] *Introduction à l'étude de la théorie des nombres et de l'algèbre supérieure*, Paris : Nony, 1895.
BOUCARD (Jenny)
[2011] Louis Poinsot et la théorie de l'ordre : un chaînon manquant entre Gauss et Galois ?, *Revue d'histoire des mathématiques*, 17, fasc. 1 (2011), p. 41-138.
BOURBAKI (Nicolas)
[1960] *Eléments d'histoire des mathématiques,* Paris : Hermann, 1960.
BRECHENMACHER (Frédéric)
[2006] *Histoire du théorème de Jordan de la décomposition matricielle (1870-1930),* Thèse de doctorat, Ecole des Hautes Etudes en Sciences sociales, Paris, 2006.
[2007] La controverse de 1874 entre Camille Jordan et Leopold Kronecker, *Revue d'Histoire des Mathématiques*, tome 13, fasc.2 (2007), p. 187-257.
 [2010] Une histoire de l'universalité des matrices mathématiques, *Revue de Synthèse,* Vol 131, n°4 (2010), p. 569-603.
[2011a] "Self-portraits with Évariste Galois (and the shadow of Camille Jordan), " *Revue d'histoire des mathématiques*, t. 17, fasc. 2, *to appear.*
[2011b] Autour de pratiques algébriques de Poincaré, *to appear.*
[2012] Galois Got his Gun, *to appear.*
BRECHENMACHER (Frédéric), EHRHARDT (Caroline)
[2010] On the identities of algebra in the 19th century, *Oberwolfach Reports*, n°12/2010, p. 24-31.
BRIOSCHI (Francesco)
[1889a] Les discriminants des résolvantes de Galois, *Comptes rendus hebdomadaires des séances de l'Académie des sciences,* t. 108 (1889), p. 878-879.

---

[50] The translation in German of Dickson's book on associative algebras in 1927 has often been considered as the dawn of a new era in the mathematical relations between the U.S.A. and Germany. For a case study of the collaboration between Albert and Hasse on division algebras, see [Fenster and Schwermer 2007].

[1889b] Sur la dernière communication d'Halphen à l'Académie, *Comptes rendus hebdomadaires des séances de l'Académie des sciences,* t. 109 (1889), p. 520-522.

BURNSIDE (William)
[1897] *Theory of Groups of Finite Order*, Cambridge: Cambridge University Press, 1897.

CLEBSCH (Alfred)
[1864] Ueber die Anwendung der Abelschen Functionen in der Geometrie, *Journal für die reine und angewandte Mathematik*, t. 63, p. 189-243.

CLEBSCH (Alfred), GORDAN (Paul)
[1866] *Theorie der Abelschen Functionen,* Leipzig : Teubner.

CORRY (Leo)
[1996] *Modern Algebra and the Rise of Mathematical Structures*, Basel: Birkhäuser, 1996.

DICKSON (Leonard E)
[1896] The analytic representation of substitutions on a power of a prime number of letters with a discussion of the linear group. Part I, II., *Annals of Math.* 11 (1896), p. 65-120, 161-183.
[1901] *Linear groups with an exposition of the Galois field theory,* Leipzig: Teubner, 1901.
[1924 /1928] A new theory of linear transformation and pairs of bilinear forms, *Proceedings Congress Toronto,* 1 (1928), p. 361-363.

DIEUDONNE (Jean)
[1962] Notes sur les travaux de Camille Jordan relatifs à l'algèbre linéaire et multilinéaire et la théorie des nombres, *in* [Jordan Œuvres, 3, p. V-XX].
[1978] *Abrégé d'histoire des mathématiques,* Paris : Hermann, 1978.

EHRHARDT (Caroline)
[2007] *Evariste Galois et la théorie des groupes. Fortune et réélaborations (1811-1910)*, Thèse de doctorat. Ecole des Hautes études en sciences sociales. Paris, 2007.

FREI (Gunther)
[2007] The unpublished section eight: On the way to function fields over a finite field, *in* [Goldstein, Schappacher, Schwermer 2007], p. 159-198.

FROBENIUS (Georg)
[1877] Ueber linear Substitutionen und bilineare Formen, *Journal für die reine und angewandte Mathematik,* t. 84 (1877), p. 1-63.
[1879] Theorie der bilinearen Formen mit ganzen Coefficienten, *Journal für die reine und angewandte Mathematik,* t. 86, p. 147-208.

FENSTER (Della-D.)
[1998] Leonard Eugene Dickson and his Work in the Arithmetics of Algebras, *Arch. Hist. Exact Sci.,* 52, p. 119-159.

FENSTER (Della-D.) and SCHWERMER (Joachim)
[2005] A Delicate Collaboration : Adrian Albert and Helmut Hasse and the Principal Theorem in Division Algebras in the Early 1930's, *Archive for History of Exact Sciences,* 59 (2005), p. 349-379.

FUCHS (Lazarus)
[1878] Ueber die linearen Differentialgleichungen zweiter Ordnung, welche algebraisch Integrale besitzen, *Mathematische Annalen*, 85 (1878), p. 1-26.

GALOIS (Évariste)
[1846] Œuvres mathématiques, *Journal de mathématiques pures et appliquées,* 11 (1846), 381–444.
[1897] *Œuvres mathématiques d'Évariste Galois, publiées sous les auspices de la Société Mathématique de France, avec une introduction par M. Émile Picard*, Paris: Gauthier-Villars, 1897.

GAUSS (Carl Friedrich)
[1801] *Disquitiones arithmeticae*, Leipzig : Fleischer, 1801.

GAUTHIER (Sébastien)
[2007] *La géométrie des nombres comme discipline*. Thèse de doctorat. Université Pierre et Marie Curie : Paris, 2007.
[2009] La géométrie dans la géométrie des nombres : histoire de discipline ou histoire de pratiques à partir des exemples de Minkowski, Mordell et Davenport, *Revue d'histoire des mathématiques,* **15-2** (2009), p. 183-230.

GIERSTER (Joseph)
[1881] Die Untergruppen der Galois'schen Gruppe der Modulargleichungen für den Fall eines primzahligen Transformationsgrades, *Mathematische Annalen*, t. 18 (1881), p. 319-365.

GISPERT (Hélène),

[1982] *Camille Jordan et les fondements de l'analyse : Comparaison de la 1ère édition (1882-1887) et de la 2ème (1893) de son cours d'analyse de l'école Polytechnique*, Thèse de doctorat, Orsay : 1982.

[1991] La France mathématique*. La Société Mathématique de France (1870-1914)*, Cahiers d'histoire et de philosophie des sciences, Paris : Belin, 1991.

GOLDSTEIN (Catherine)

[1995] *Un théorème de Fermat et ses lecteurs*, Saint-Denis : PUV, 1995.

[1999] Sur la question des méthodes quantitatives en histoire des mathématiques : le cas de la théorie des nombres en France (1870- 1914), *Acta historiae rerum necnon technicarum*, nouv. sér., vol. 3 (1999), p. 187-214.

[2007] The Hermitian Form of Reading the *Disquisitiones, in* [Goldstein, Schappacher, Schwermer 2007], p. 377-410.

[2009] La théorie des nombres en France dans l'entre-deux-guerres : De quelques effets de la première guerre mondiale, *Revue d'histoire des sciences,* 62-1 (2009), p. 143-175.

[2011] Hermite's strolls in Galois fields, *Revue d'histoire des mathématiques,* to appear.

GOLDSTEIN (Catherine), GRAY, (Jeremy), RITTER (Jim) (eds.)

[1996] *L'Europe mathématique. Histoires, mythes, identités*, Paris : Éditions de la MSH, 1996.

GOLDSTEIN (Catherine), SCHAPPACHER (Norbert)

[2007a] A Book in Search of a Discipline (1801-1860) *in* [Goldstein, Schappacher, Schwermer 2007], p. 3-66.

[2007b] Several Disciplines and a Book (1860–1901), *in* [Goldstein, Schappacher, Schwermer 2007], p. 67-104.

GOLDSTEIN (Catherine), SCHAPPACHER (Norbert), SCHWERMER (Joaquim) (eds.)

[2007] *The Shaping of Arithmetics after C. F. Gauss's Disquisitiones Arithmeticae*, Berlin: Springer, 2007.

GRAY (Jeremy)

[2000] *Linear differential equations and group theory from Riemann to Poincaré,* 2$^{de}$ ed. Boston: Birkhäuser, 2000.

GUGGENHEIMER (Heinrich)

[1977] The Jordan curve Theorem and an Unpublished Manuscript by Max Dehn," *Archive for History of Exact Sciences,* 17 (1977), 193-200.

HALPHEN (Georges)

[1889] Sur la résolvante de Galois dans la division des périodes elliptiques par 7, *Comptes rendus hebdomadaires des séances de l'Académie des sciences,* t. 108 (1889), p. 476-477.

HAWKINS (Thomas)

[1977] Weierstrass and the Theory of Matrices, *Archive for History of Exact Sciences,* 17 (1977), p. 119-163.

[2008] Frobenius and the symbolical algebra of matrices, *Archive for History of Exact Sciences,* 62 (2008), p. 23-57.

HILBERT (David)

[1894] Grundzüge einer Theorie des Galois'schen Zahlkörpers, *Göttingen Nachrichten*, 1894, p. 224-236.

[1897] Die Theorie der Algebraischen Zahlkörper*, Jahresbericht der Deutschen Mathematiker Vereinigung*, vol. 4 (1897), p. 175-546.

[1902] Ueber die Grundlagen der Geometrie, *Göttingen Nachrichten*, 1902, p. 233-241.

[1903] Ueber die Grundlagen der Geometrie, *Mathematische annalen*, 56 (1903), p. 381-422.

HÖLDER (Otto)

[1889] Zurückführung einer beliebigen algebraischen Gleichung auf eine Kette von Gleichungen, *Mathematische Annalen*, vol. 34 (1889), p. 26-56.

JORDAN (Camille)

[1860] *Sur le nombre des valeurs des fonctions*, Thèses présentées à la Faculté des sciences de Paris par Camille Jordan, 1re thèse, Paris : Mallet-Bachelier, 1860.

[1866] Note sur les irrationnelles algébriques, *Comptes rendus hebdomadaires des séances de l'Académie des sciences,* t. 63 (1866), p. 1063-1064.

[1868a] Sur la résolution algébrique des équations primitives de degré $p^2$, *Journal de mathématiques pures et appliquées,* 32 (2) (1868), p. 111-135.

[1868b] Note sur les équations modulaires, *Comptes rendus hebdomadaires des séances de l'Académie des sciences,* t.66 (1868), p. 308-312.

[1869a] Commentaire sur Galois, *Mathematische Annalen*, vol. 1 (1869), p. 141-160.

[1869b] Sur les équations de la géométrie, *Comptes rendus hebdomadaires des séances de l'Académie des sciences*, t.68 (1869), p. 656-659.

[1869c] Sur la trisection des fonctions abéliennes et sur les vingt-sept droites des surfaces du troisième ordre, *Comptes rendus hebdomadaires des séances de l'Académie des sciences*, t.68 (1869), p. 865-869.

[1870] *Traité des substitutions et des équations algébriques,* Paris, 1870.

[1874] Mémoire sur une application de la théorie des substitutions à l'étude des équations différentielles linéaires, *Bulletin de la société mathématique de France,* 2 (1874), p. 1000-1027. [1878]

[1876] Sur les équations linéaires du second ordre dont les intégrales sont algébriques, *Comptes rendus hebdomadaires des séances de l'Académie des sciences,* t. 82 (1876), p. 605-607.

[1878] Mémoire sur les équations différentielles linéaires à intégrale algébrique, *Journal für die reine und angewandte Mathematik,* t. 84 (1878), p. 89-215.

[1880] Mémoire sur l'équivalence des formes, *Journal de l'École polytechnique,* 28 (1880), p. 111-150.

[1961-1964] *Œuvres de Camille Jordan.* Publiées sous la direction de M. Gaston Julia, par M. Jean Dieudonné. Paris : Gauthier-Villars, 1961-1964.

KIERNAN (Melvin)

[1971] The Development of Galois Theory from Lagrange to Artin, *Archive for History of Exact Sciences*, vol. 8, n° 1-2 (1971), p. 40-152

KLEIN (Felix)

[1879a] Ueber die Transformation der elliptischen Functionen und die Auflösung der Gleichungen fünften Grades, *Mathematische Annalen,* 14 (1879), p. 111-172.

[1879b] Ueber die Erniedrigung der Modulargleichungen, *Mathematische Annalen,* 14 (1879), p. 417-427.

[1884] *Vorlesungen über das Ikosaeder und die Auflösung der Gleichungen vom fünften Grade* , Leipzig, Teubner, 1884.

[1893] *Einleitung in die höhere Geometrie, I und II. Vorlesung gehalten im W. S. 1892/3 und S. S. 1893. Ausgearbeitet von F. Schilling.*Göttingen. Authogr, 1893

[1921-1923] *Gesammelte mathematische Abhandlungen,* Springer: Berlin, 1921.3. Vol.

KRONECKER (Leopold)

[1882] Grundzüge einer arithmetischen Theorie der algebraischen Grössen, *Journal für die reine und angewandte Mathematik*, t. 92 (1882), p. 1-122.

LEBESGUES (Henri).

[1923] Notices d'histoire des mathématiques. Notice sur la vie et les travaux de Camille Jordan, *L'enseignement mathématique* (1923)*,* p. 40-49.

LIE (Sophus)

[1895] Influence de Galois sur le développement des mathématiques, in DUPUY (Paul) (ed.), *Le Centenaire de l'École Normale 1795-1895*, Paris, Hachette, 1895.

MARIE (Joseph)

[1873] Des residus relatifs aux asymptites. Classification des courbes algébriques, *Comptes renuds de l'Académie des sciences de Paris*, t. 76 (1873), p. 943-947.

MARTINI (Laura),

[1999] The First Lectures in Italy on Galois Theory : Bologna, 1886- 1887, *Historia Mathematica*, vol. 26, n° 3, 1999, p. 201-223.

MINKOWSKI (Hermann)

 [1885] Untersuchungen über quadratische Formen. 1. Bestimmung der Anzahl verschiedener Formen, welche ein gegebenes Genus enthält, *Acta Mathematica*, 7 (1885), p. 201-258.

MOORE (Eliakim H.)

[1893] A doubly-infinite system of simple groups, *Bulletin of the New York Mathematical Society,* III (3) (1893), p. 73-78.

NETTO (Eugen)

[1882] *Subtitutionentheorie und ihre Anwendung auf die Algebra*, Leipzig, Teubner, 1882.

NEUMANN (Olaf)

[1997] Die Entwicklung der Galois-Theorie zwischen Arithmetik und Topologie (1850 bis 1960). *Archive for History of Exact Sciences,* 50 (1997), 291–329.

[2007] The *Disquisitiones Arithmeticae* and the Theory of Equations, in [Goldstein, Schappacher, Schwermer 2007], p. 107-128.

NEUMANN (Peter M.)

[2006] The concept of Primitivity in Group Theory and the Second Memoir of Galois, *Archive for History of Exact Sciences,* 60 (2006), p. 379-429.

NICHOLSON (Julia)
[1993] The development and understanding of the concept of quotient group, *Historia Mathematica*, 20 (1993), p. 68-88.

PARSHALL (Karen H.)
[1985] J. H.M. Wedderburn and the Structure Theory of Algebras, *Archive for History of Exact Sciences*, 32 (1985), p. 223–349.
[2004] Defining a mathematical research school: the case of algebra at the University of Chicago, 1892-1945, *Historia Mathematica*, 31 (2004), p. 263-278.

PARSHALL (Karen H.), ROWE (David E.)
[1994] *The Emergence of the American Mathematical Research Community (1876–1900)*: *J.J. Sylvester, Felix Klein, and E.H. Moore.* Providence, American Mathematical Society, London: AMS/LMS Series in the History of Mathematics, vol. 8, p. 1994.

PETRI (Birgit), SCHAPPACHER (Norbert)
[2004] From Abel to Kronecker. Episodes from 19th Century Algebra, p. 227-266 in LAUDAL (Olav Arnfinn), PIENE,(Ragni), *The Legacy of Niels Henrik Abel ; the Abel Bicentennial, Oslo 2002*, Springer Verlag, Berlin, 2004.
[2007] On Arithmetization, in [Goldstein, Schappacher, Schwermer 2007, p. 343-374].

PICARD (Émile)
[1896] *Traité d'analyse*, t. III, Paris : Gauthier-Villars, 1896.
[1922] Résumé des travaux mathématiques de Jordan, *Comptes rendus de lAacadémie des sciences de Paris,* t. 174 (1922), p. 210-211.

PIERPONT (James)
[1897] Early History of Galois Theory of Equations, *Bulletin of the American Mathematical Society*, vol. 2, n° 4 (1897), p. 332-340.

POINCARE (Henri)
[1883] Sur l'intégration algébrique des équations linéaires, *Comptes rendus hebdomadaires des séances de l'Académie des sciences,* t. 97 (1883), p. 984-985.

POINSOT (Louis)
[1808] Commentaire sur le livre de Lagrange, *Magasin encyclopédique*, t. 4, juill.-août 1808, p. 343-375.
[1845] Réflexion sur les principes fondamentaux de la théorie des nombres, *Journal de Mathématiques Pures et Appliquées,* 10 (1845), p. 1–101.
[1851] Théorie Nouvelle de la rotation des corps, *Journal de mathématiques pures et appliquées*, XVI (1851), p. 73-129, 289-336.

ROWE (David E.)
[1989] The early geometrical works of Sophus Lie and Felix Klein, *The history of modern mathematics* (Boston, MA, 1989), p. 209-273.

SCHÖNFLIES (Arthur)
[1891a] Sur les équations de deux surfaces minima périodiques, possédant la symétrie de l'octaèdre, *Comptes rendus de l'Académie des sciences de Paris*, 112 ( 1891), p. 515-518.
[1891b] *Krystallsysteme und Krystallstructur,* Leipzig : Teubner, 1891.
[1892] Ueber gewisse geradlinig bergenzte Stücke Riemann'scher Flächen, *Göttingen Nachrichten* (1892), p. 257-267.
[1893] *La géométrie du mouvement. Traduit de l'allemand par Ch. Speckel,* Paris : Gauthier-Villard, 1893.
[1896] Ueber einen Satz aus der Analysis situs, *Göttingen Nachrichten* (1896), p. 79-89

SCHOLZ (Erhard)
[1989] *Symmetrie, Gruppe, Dualität. Zur Beziehung zwischen theoretischer Mathematik und Anwendungen in Kristallographie und Baustatik des 19. Jahrhunderts*, Basel : Birkhäuser, 1989.

SCHWARZ (Hermann)
[1872] Ueber diejenigen Fälle, in welchen die Gaussische hypergeometrische Reihe eine algebraische Function ihres vierten Elementes darstellt, *Journal für die reine und angewandte Mathematik,* 75 (1873), p. 292-335.

SÉGUIER (Jean-Armand de)
[1904a] *Théorie des groupes finis. Éléments de la théorie des groupes abstraits*, Paris : Gauthier-Villars, 1904.

[1907] Sur la théorie des matrices, *Comptes rendus hebdomadaires des séances de l'Académie des sciences,* t. 145 (1907), p. 1259-1260.

SERRET (Joseph Alfred)

[1866] *Cours d'algèbre supérieure*, 3e éd., Paris : Gauthier-Villars, 1864, 2 vols.

SIEGMUND-SCHULTZE (Reinhard)

[1993]. *Mathematische Berichterstattung in Hitler- deutschland. Der Niedergang des "Jahrbuchs über die Fortschritte der Mathematik"*, Studien zur Wissenschafts-, Sozial- und Bildungsgeschichte der Mathematik 9, Vandenhoeck & Ruprecht, Göttingen.

TATON (René)

[1947] Les relations scientifiques d'Évariste Galois avec les mathématiciens de son temps, *Revue d'histoire des science*s, 1 (1947), p. 14–130.

TAZZIOLI (Rossana)

[1994] Schwarz's critique and interpretation of the Riemann representation theorem, *Rendiconti del circolo matematico del Palermo*, (2) 34 (1994), 95-132.

TIGNOL (Jean-Pierre)

[2001] *Galois' Theory of Algebraic Equations*, World Scientific Publishing Company, 2001.

VILLAT (Henri)

[1922] Camille Jordan, *Journal de mathématiques pures et appliquées*, (9) 1 (1922), p. 1-5.

WAERDEN (Bartel, van der)

[1985] *A history of Algebra : from Al-Khwàrizmi to Emmy Noether*, New York, Springer Verlag, 1985

WEBER (Heinrich)

[1893] Die allgemeinen Grundlagen der Galois'schen Gleichungstheorie, *Mathematische Annalen*, vol. 43 (1893), p. 521-549.

[1895] *Lehrbuch der Algebra*, Braunschweig, F. Vieweg und Sohn, 1895-1896, 2 vols.

WUSSING (Hans)

[1984] *The Genesis of the Abstract Group Concept*, MIT Press, Cambridge: Mass.,1984.

**Annex**

**I. *Jahrbuch* classifications**

The *Jahrbuch* mathematical classification changes in time.

It is however not compulsory to follow closely the evolution of the *Jahrbuch* for investigating the Jordan's works and their reception. Indeed, only a few sections of the *Jahrbuch* classifications played an important role in regard with Jordan's works.

I am therefore simplifying here the evolution of the *Jahrbuch* classifications. First, I am keeping the same numeration when a section evolves in continuity (for instance both II. Algebra 3. Substitutions and III. Algebra. 3. Group theory will be referred to as II.3.).

Section that did not exist at the beginning of the time period are referred to by adding a decade to their original numeration. For instance, prior to 1914, the third section of the *Jahrbuch* classification was referring to number theory. After 1914, the third section was referring to set theory. In the present paper, I am therefore referring to set theory as section 13.

**I. History and philosophy.**

**II. Algebra**. 1. Equations. 2. Substitutions ( turns into "group theory after 1895). 3. Theory of forms.

**III. Number theory**. 1. Generalities. 2. Theory of forms (Arithmetic and Algebra after 1914).

**IV. Series**.

**V. Probability**.

**VI. Differential and integral calculus**. 3. Integral calculus. 4. Definite Integral. 5. Ordinary differential equations.

**VII. Function theory**. 1. Generalities. 2. Special functions.

**VIII. Pure, elementary, and synthetic geometry**. 1. Principles of geometry. 2.  Continuity, analysis situs. 5. New analytic geometry.

**IX. Analytic geometry**. 3. Analytic geometry of space. 5. Correspondences, transformations.XI.

**X. Mechanic**. 2. Cinematic.

**XI. Mathematical physics**. 1. Molecular physics.

**XIII. Set theory** (after 1914)

**XIV. Analysis** (after 1914). 3. Real function theory. 4. Complex function theory. 5. Conformal mapping. 8. Continuous groups. 9. Differential equations. 11. Differential calculus. 13. Potential theory, partial differential equations. 15. Variations calculus.

**XV. Geometry** (after 14). 1. Foundations of geometry, 2. Topology, 5. Analytic geometry, 6. Differential geometry.

**XXIII. Arithmetic and algebra**. 2. Determinant and matrices.

## II. Distributions of the *Jahrbuch* classifications of Jordan's works

Diagrams 1.

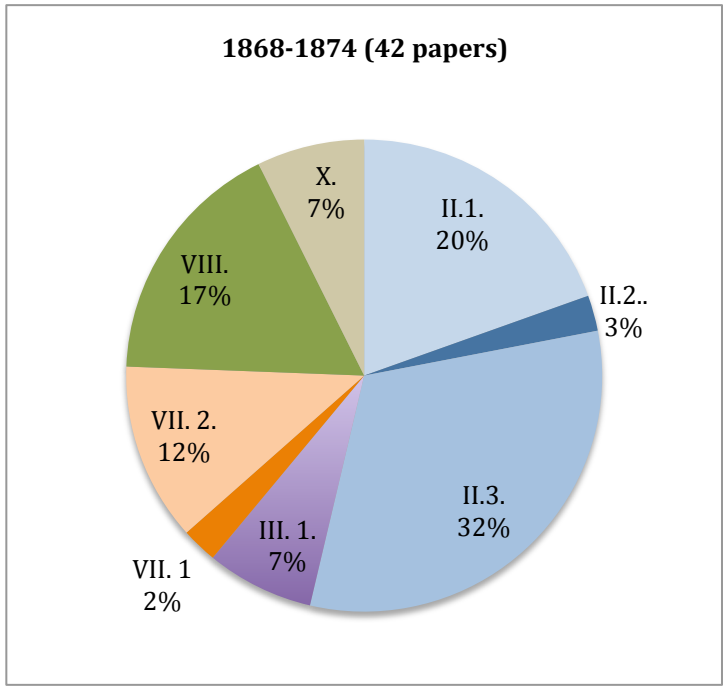

**1868-1874 (42 papers)**

X. 7%
II.1. 20%
II.2.. 3%
II.3. 32%
III. 1. 7%
VII. 1 2%
VII. 2. 12%
VIII. 17%

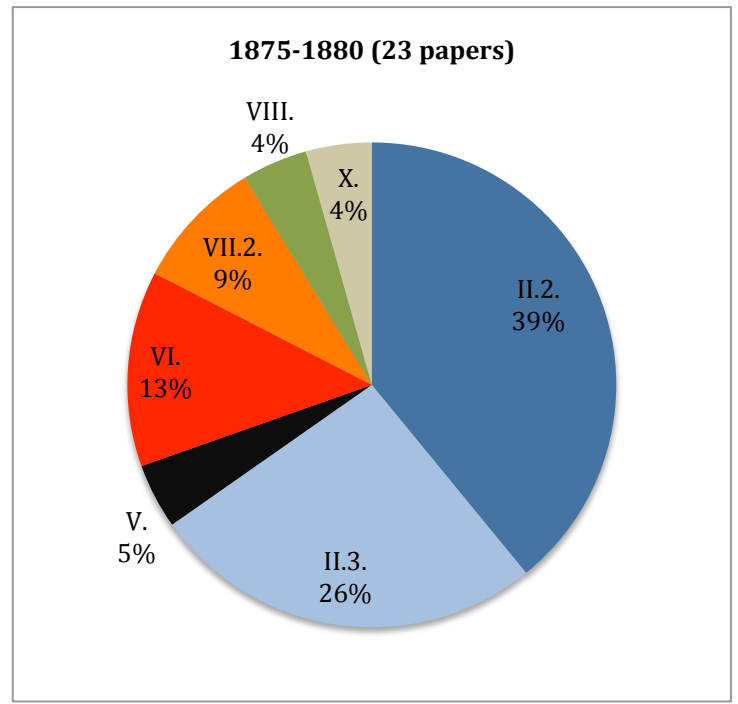

**1875-1880 (23 papers)**

VIII. 4%
X. 4%
VII.2. 9%
VI. 13%
II.2. 39%
V. 5%
II.3. 26%

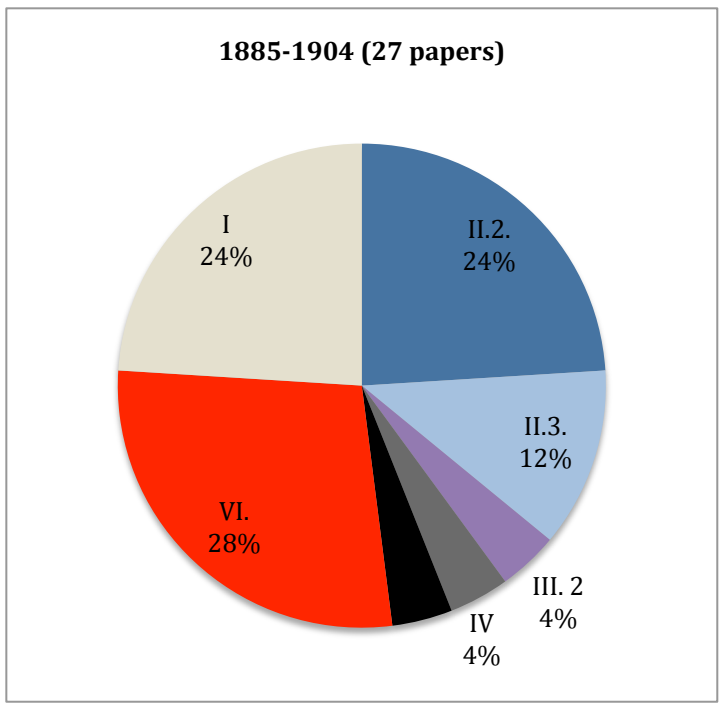

**1885-1904 (27 papers)**

I 24%
II.2. 24%
II.3. 12%
III. 2 4%
IV 4%
VI. 28%

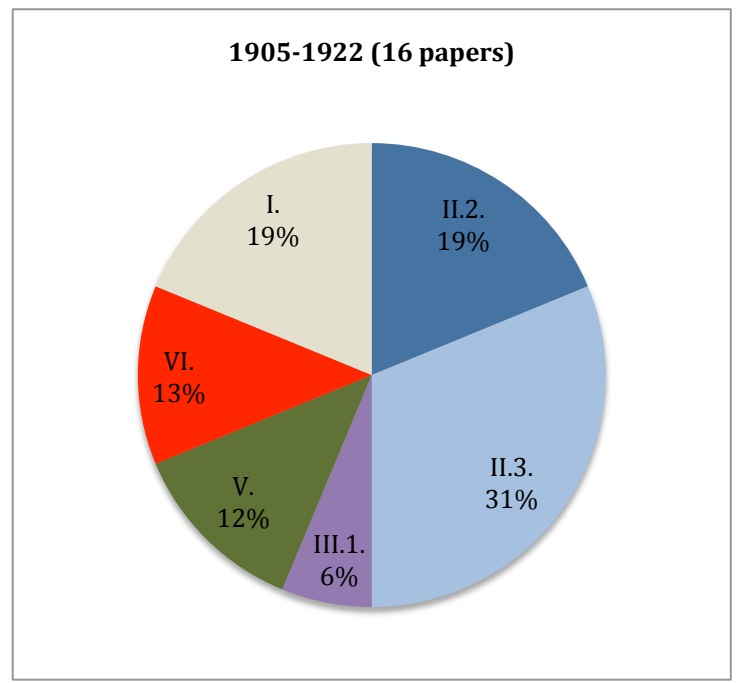

**1905-1922 (16 papers)**

I. 19%
II.2. 19%
VI. 13%
II.3. 31%
V. 12%
III.1. 6%

## II. The Jordan corpus
## II. 1. An overview on the Jordan corpus

Diagrams n° 2.

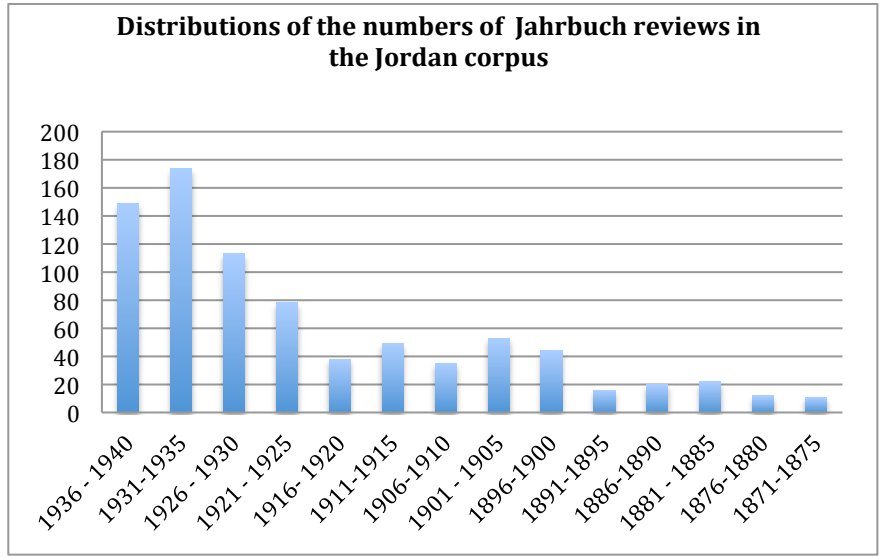

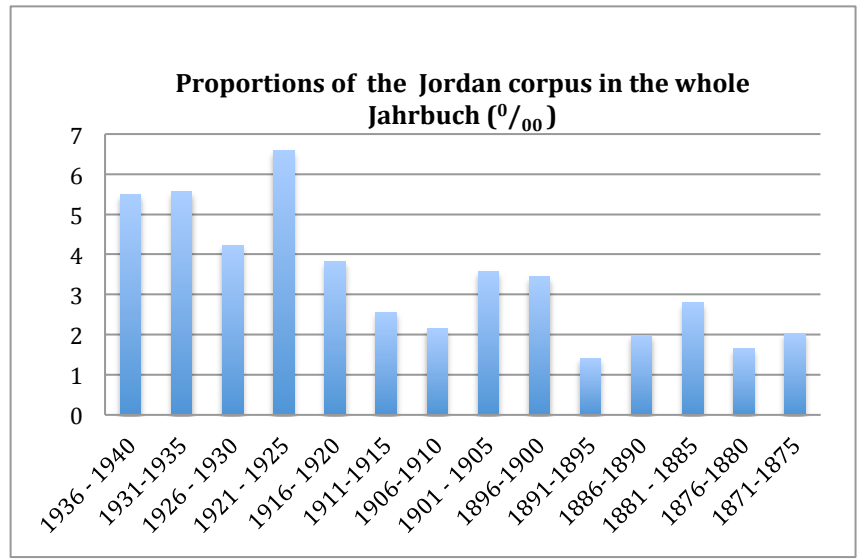

## II.2. Distributions of the classifications in the Jordan corpus
### Diagrams n°3

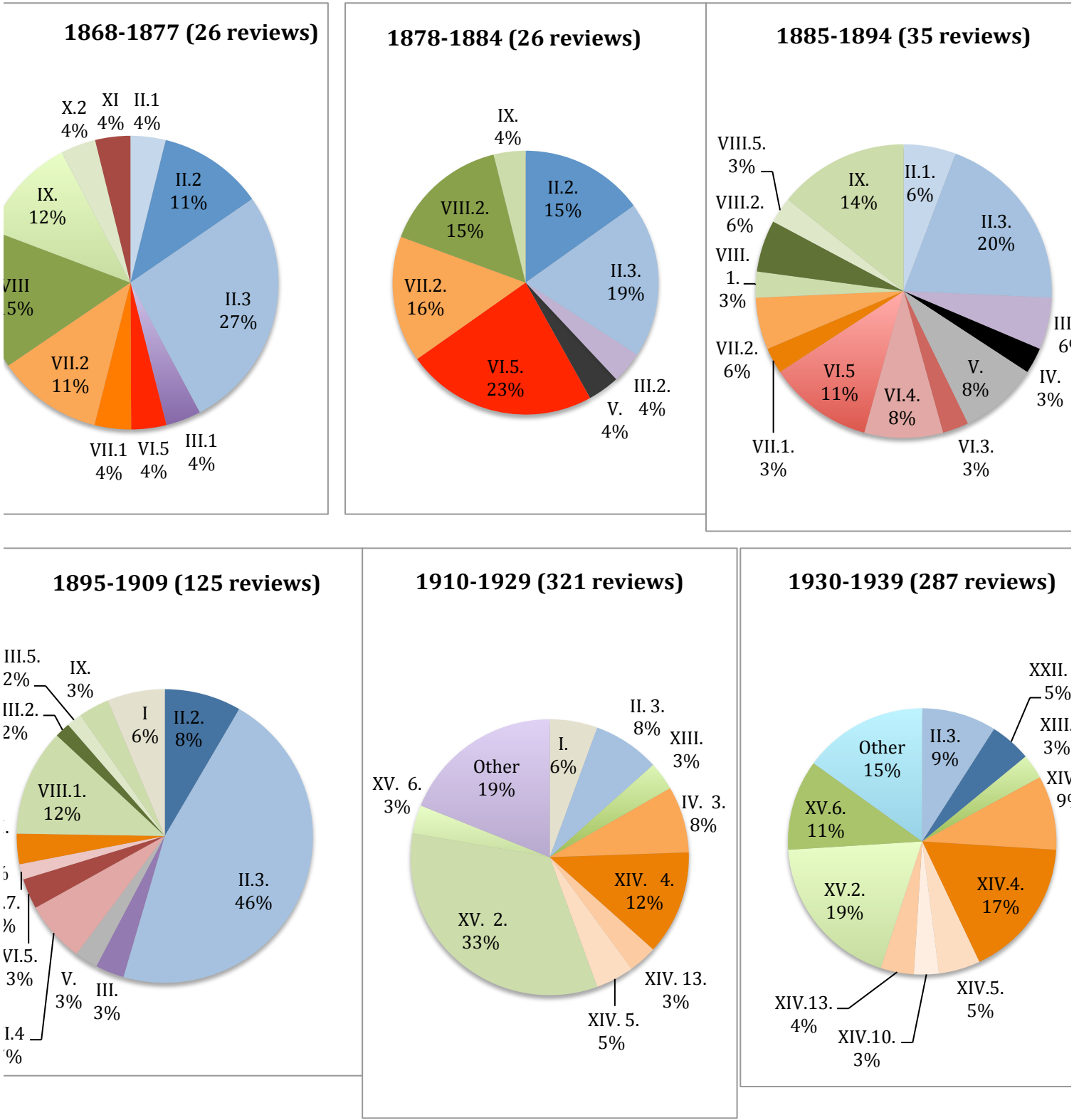

## II.3. Distibutions of the languages in the Jordan corpus
### Diagrams n°4

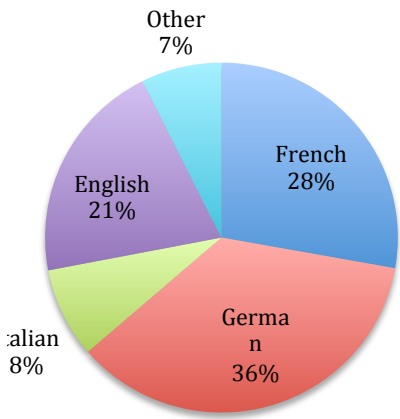

**1868-1877**

Other 7%
French 28%
English 21%
Italian 8%
German 36%

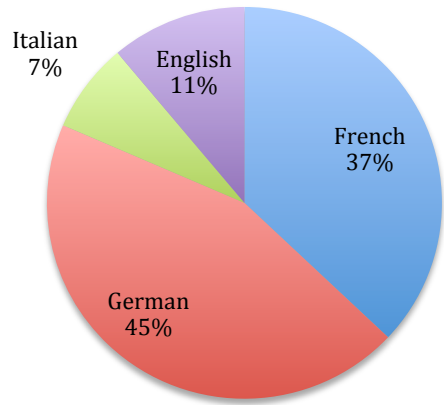

**1878-1884**

Italian 7%
English 11%
French 37%
German 45%

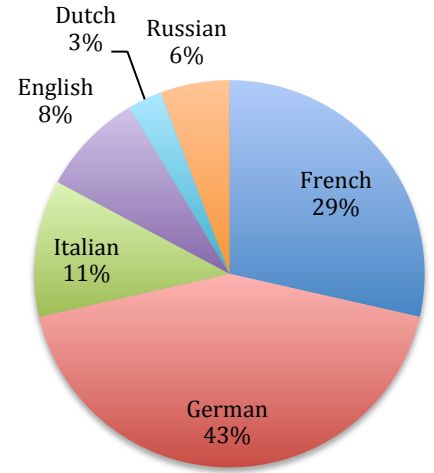

**1885-1894**

Dutch 3%
Russian 6%
English 8%
French 29%
Italian 11%
German 43%

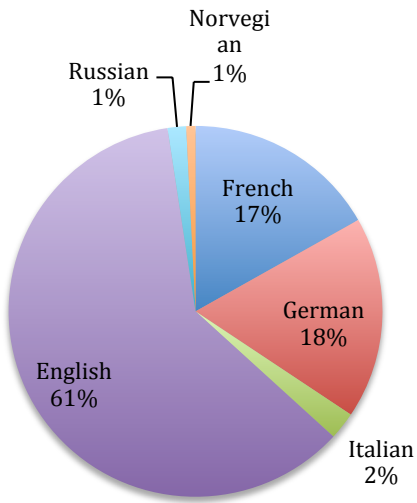

**1895-1909**

Norvegian 1%
Russian 1%
French 17%
German 18%
English 61%
Italian 2%

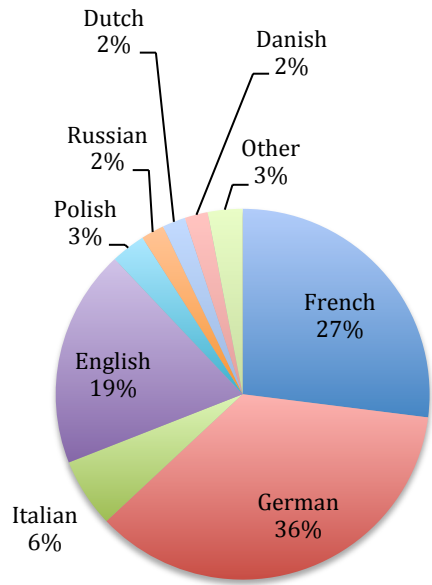

**1910-1929**

Dutch 2%
Danish 2%
Russian 2%
Other 3%
Polish 3%
French 27%
English 19%
Italian 6%
German 36%

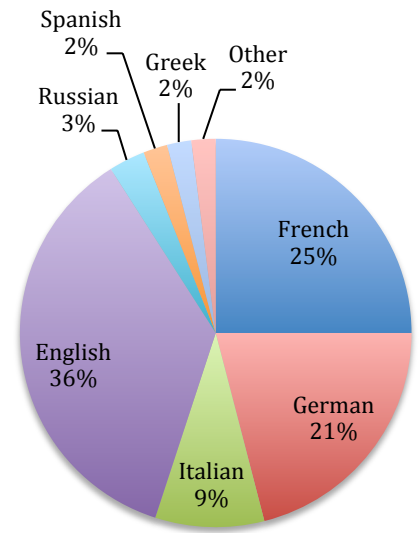

**1930-1939**

Spanish 2%
Greek 2%
Other 2%
Russian 3%
French 25%
English 36%
German 21%
Italian 9%

## II.3. Types of references to Jordan
### Diagrams n°5

### 1868-1877

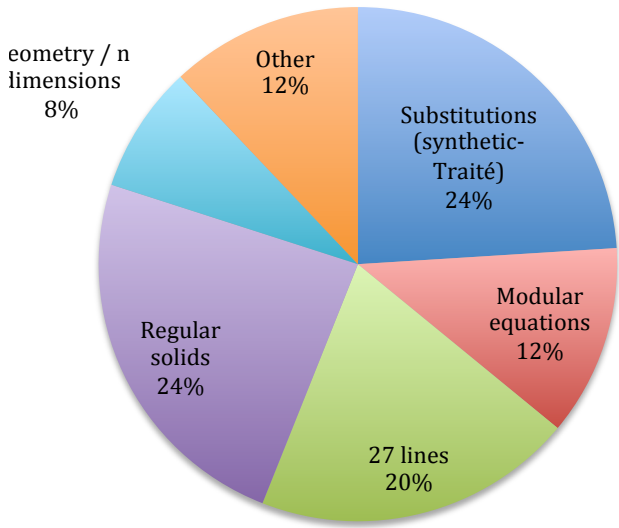

- Substitutions (synthetic-Traité) 24%
- Modular equations 12%
- 27 lines 20%
- Regular solids 24%
- Geometry / n dimensions 8%
- Other 12%

### 1878-1884

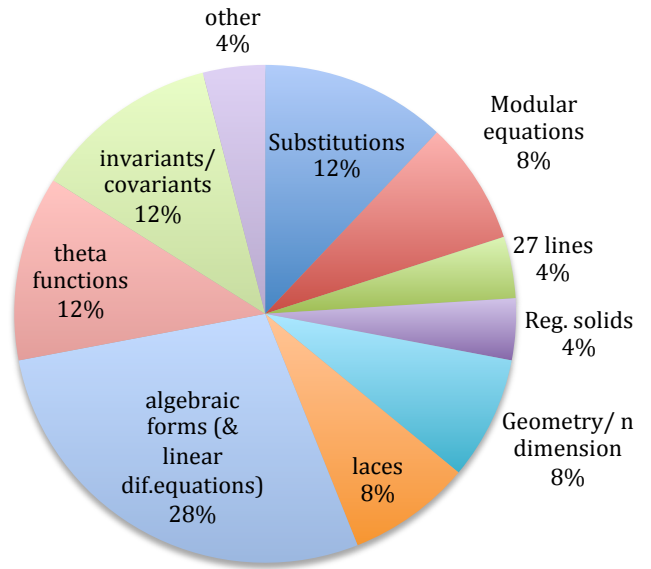

- other 4%
- Substitutions 12%
- Modular equations 8%
- 27 lines 4%
- Reg. solids 4%
- Geometry/ n dimension 8%
- laces 8%
- algebraic forms (& linear dif.equations) 28%
- theta functions 12%
- invariants/covariants 12%

### 1885-1894

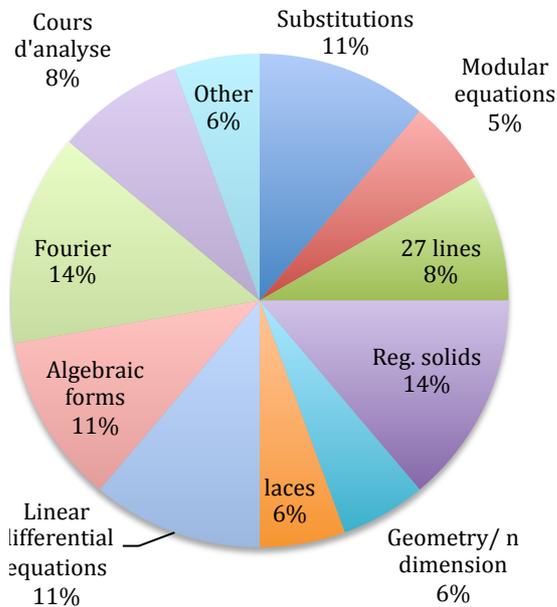

- Cours d'analyse 8%
- Substitutions 11%
- Modular equations 5%
- Other 6%
- 27 lines 8%
- Reg. solids 14%
- Geometry/ n dimension 6%
- laces 6%
- Linear differential equations 11%
- Algebraic forms 11%
- Fourier 14%

### 1895-1909

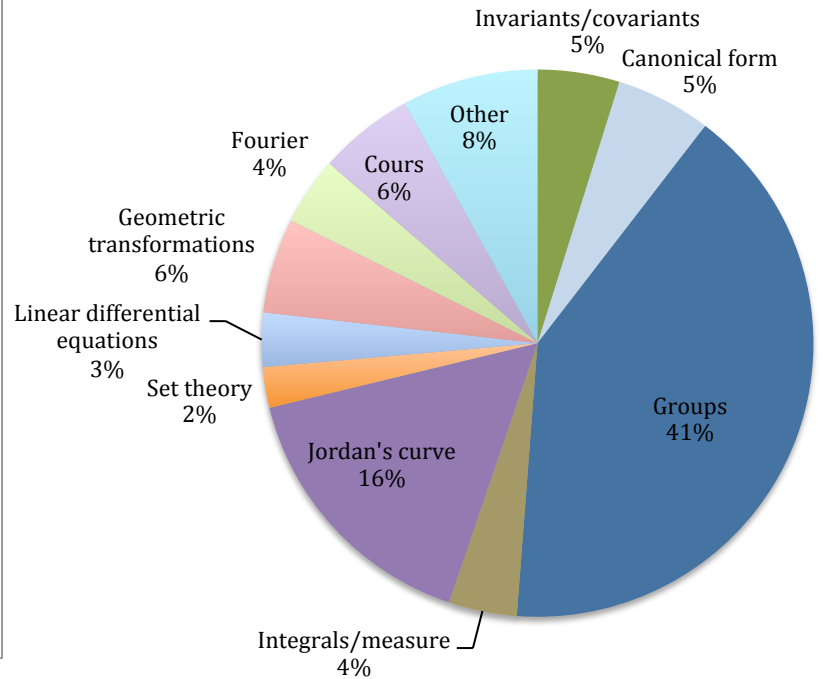

- Invariants/covariants 5%
- Canonical form 5%
- Fourier 4%
- Cours 6%
- Other 8%
- Geometric transformations 6%
- Linear differential equations 3%
- Set theory 2%
- Jordan's curve 16%
- Groups 41%
- Integrals/measure 4%

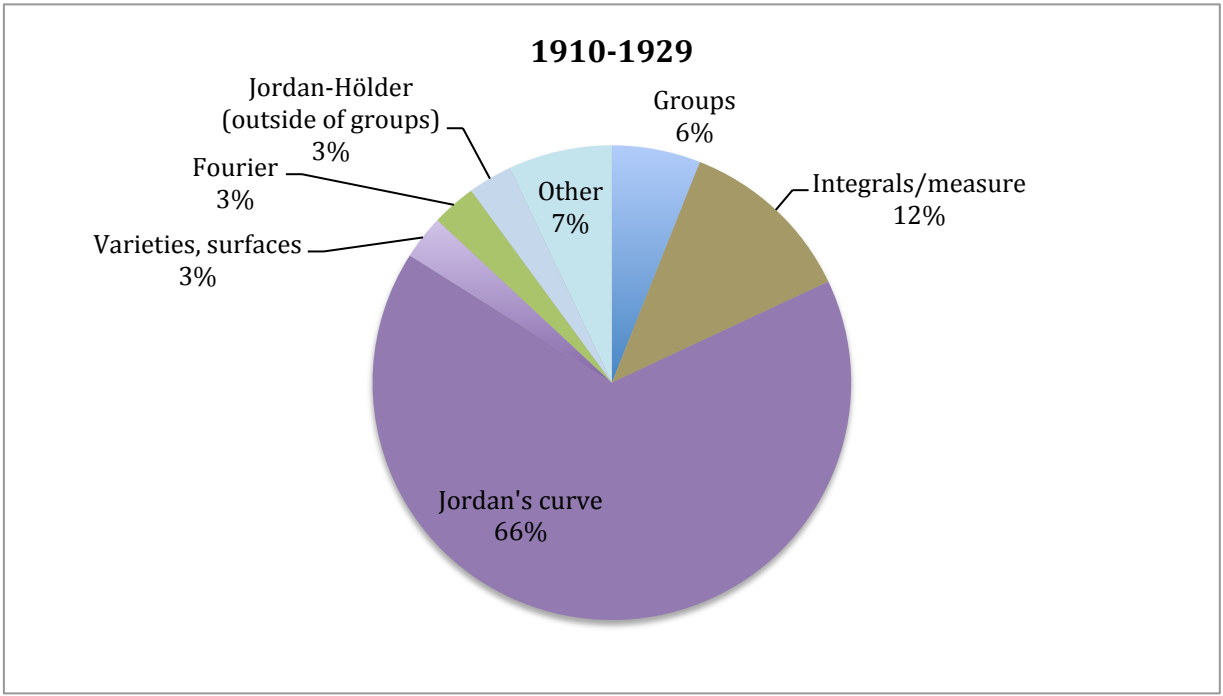

**1910-1929**

Jordan-Hölder (outside of groups) 3%
Fourier 3%
Varieties, surfaces 3%
Other 7%
Groups 6%
Integrals/measure 12%
Jordan's curve 66%

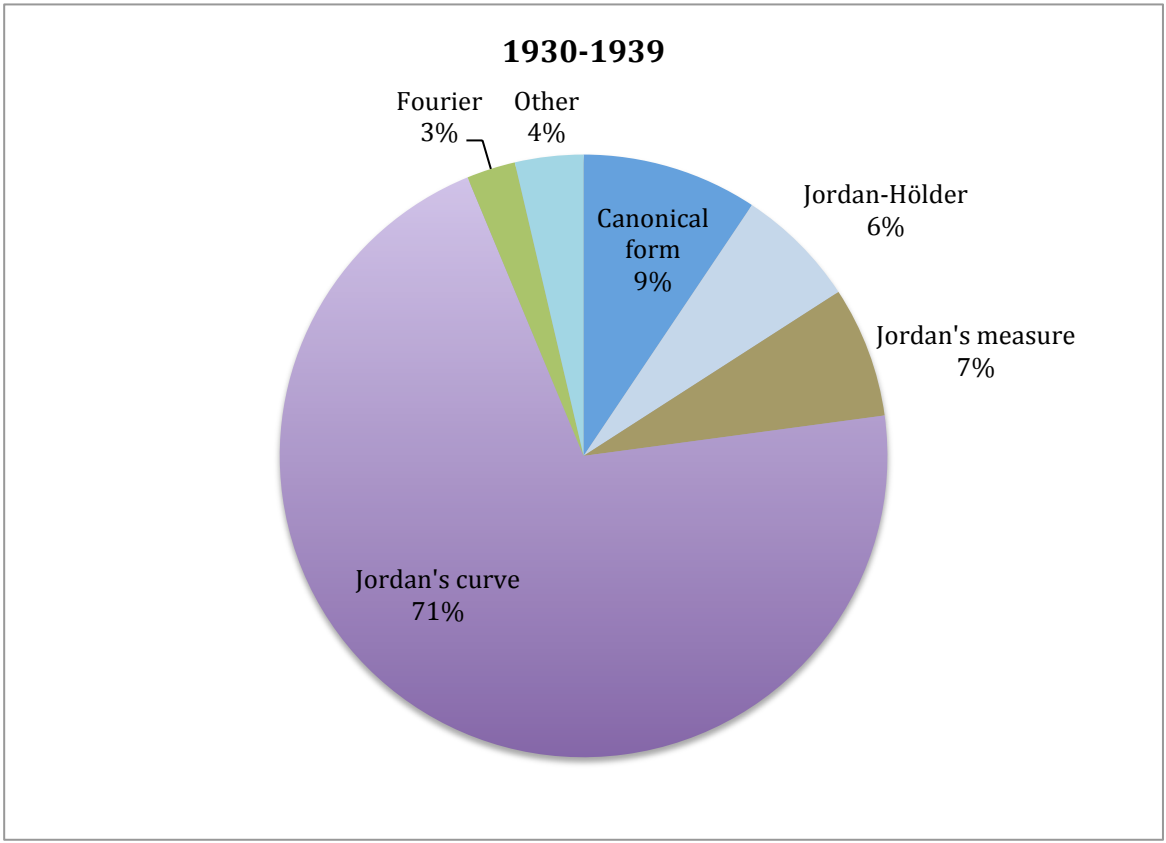

**1930-1939**

Fourier 3%
Other 4%
Canonical form 9%
Jordan-Hölder 6%
Jordan's measure 7%
Jordan's curve 71%

# III. The Galois corpus

## III. 1 An overview (proportions of the Galois corpus in the *Jahrbuch* ($^0/_{00}$)

Diagram n°6

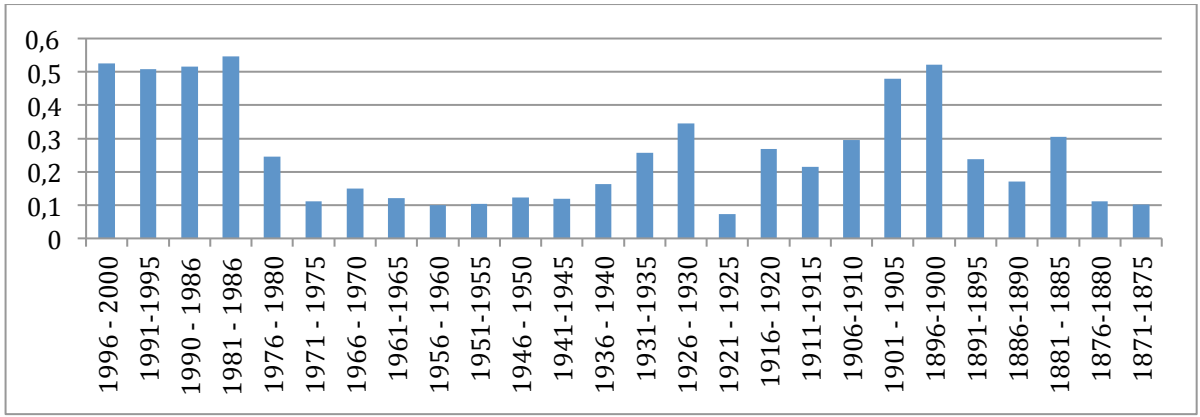

## III. 2 Distributions of the classifications in the Galois corpus

Diagrams n° 7

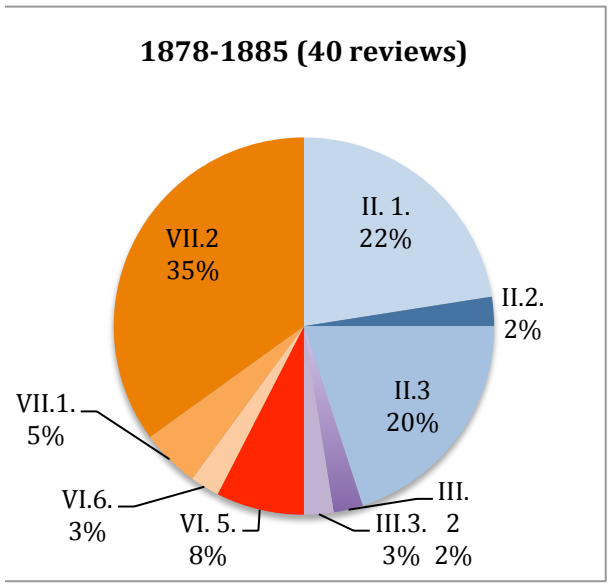

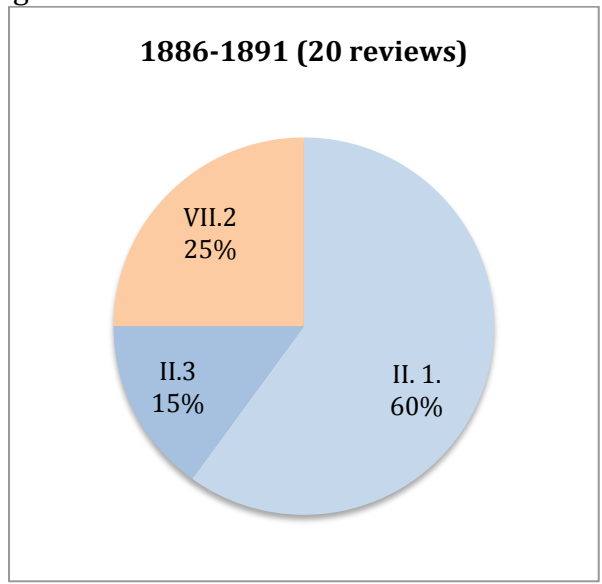

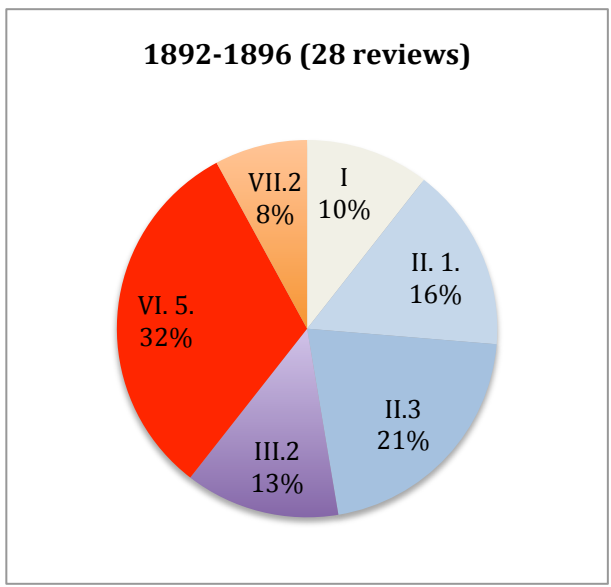

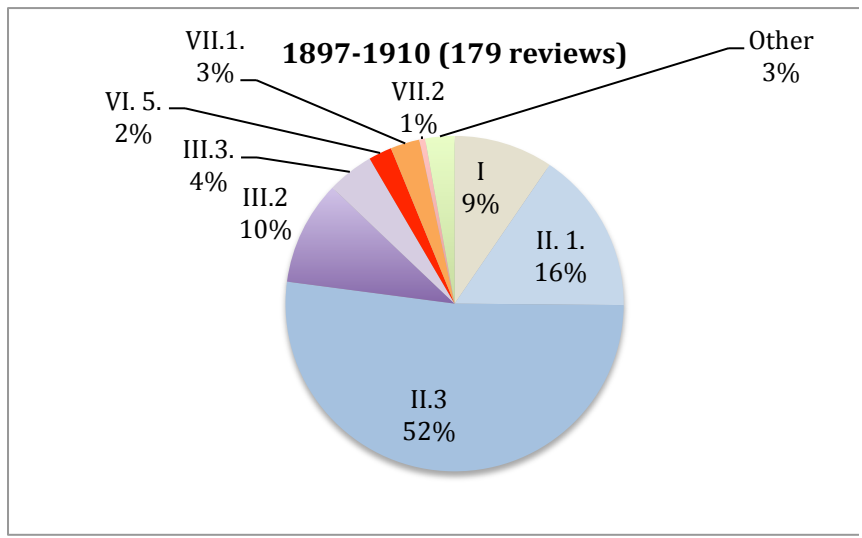